\documentclass[12pt]{amsart}
\usepackage{epsfig}

 \newcommand{\new}{\newcommand}                        
 \new{\trunc}{{\otimes}}                               
 \new{\tnsr}{\otimes}                                  
 \new{\tensor}{\otimes}                                
 \new{\iso}{\cong}                                     
 \new{\union}{\cup}                                    
 \new{\W}{\mathfrak{W}}                                
 \new{\g}{\mathfrak{g}}                                
 \new{\h}{\mathfrak{h}}                                
 \new{\uqg}{U_q(\g)}                                   
 \new{\bracket}[1]{\langle#1\rangle}                   
 \new{\qdim}{\operatorname{qdim}}                      
 \new{\Hom}{\operatorname{Hom}}                        
 \new{\C}{\Bbb{C}}                                     
 \new{\R}{\Bbb{R}}                                     
 \new{\Cat}{\mathfrak{C}}                              
 \new{\ZZ}{\mathfrak{Z}}                               
 \new{\F}{\mathfrak{F}}                                
 \new{\Z}{\Bbb{Z}}                                     
 \new{\Q}{\Bbb{Q}}                                     
 \new{\N}{\Bbb{N}}                                     
 \newcounter{letter}                                   
 \newtheorem{thm}{Theorem}                             
 \newtheorem{theorem/definition}{Theorem/Definition}   
 \newtheorem{prop}{Proposition}                        
 \newtheorem{cor}{Corollary}                           
 \newtheorem{lem}{Lemma}                               
\theoremstyle{definition}                              
\newtheorem{rem}{Remark}                               

\new{\pic}[5]{\raisebox{#3pt}{
\hspace{#4pt}\epsfig{file=#1.ps,height=#2pt}\hspace{#5pt}}}
\begin{document}
\title[Nonsimply Connected Lie Groups]{Jones-Witten  Invariants for
Nonsimply-Connected Lie Groups and The Geometry of the Weyl Alcove}
\author{Stephen F. Sawin}
\address{Fairfield University\\(203)254-4000x2573\\
ssawin@fair1.fairfield.edu}
\maketitle

\begin{abstract}
The quotient process of M\"uger and Brugui\`eres is used to construct modular categories
and TQFTs out of closed subsets of the Weyl alcove of a simple Lie
algebra. In particular it is determined at 
which levels closed subsets associated to nonsimply-connected groups
lead to TQFTs.   Many of
these TQFTs are shown to decompose into a tensor product of TQFTs coming from smaller
subsets.  The ``prime" subsets among these are  classified, and apart from some giving
TQFTs depending on homology as described by Murakami, Ohtsuki and Okada, they are shown to
be in one-to-one correspondence with the TQFTs predicted by Dijkgraaf and Witten to be
associated to Chern-Simons theory with a nonsimply-connected Lie group.  Thus in
particular a rigorous construction of the Dijkgraaf-Witten TQFTs is
given.  As a byproduct, a purely 
quantum groups proof of the modularity of the full Weyl alcove for arbitrary quantum
groups at arbitrary levels is given.
\end{abstract}

\section*{Introduction}

Since Witten's seminal paper \cite{Witten89a} relating the Jones
polynomial \cite{Jones85} to Chern-Simons field theory, the link and
three-manifold invariants descendent from the Jones polynomial have admitted
two apparently incompatible interpretations.  On the one hand all can be
defined combinatorially in terms of quantum groups
\cite{Kirillov96,RT90,RT91}, an algebraic language for rigorously and
coherently computing them and proving their basic properties.   Unfortunately,
in this framework it is very difficult to relate the invariants to classical
topology and geometry and, partly as a consequence, the invariants have
answered very few questions which might have been asked before their
appearance, the ultimate test of the significance of a new field.  On the
other hand they can be described geometrically as an ill-defined average over
the space of connections \cite{Witten89a}.   This definition offers a
beautiful and compelling intrinsically three-dimensional framework for the
invariants which connects them to much of the exciting geometry and physics
that has arisen over the past few decades.  But this definition is
completely nonrigorous  
because of its reliance on the path integral, a heuristic technique of physics
whose precise mathematical formulation is widely believed to be a
problem we will leave to 
our grandchildren.

 Perhaps the central problem of the subject is to unite these two viewpoints.  
A complete resolution of this problem would amount to a rigorous
interpretation of the path integral in this particular case, and while this is
arguably easier than such an interpretation in more general or more physically
interesting situations, it  should be viewed as a very long term goal.  
Still, much interesting progress has been made towards the goal of putting
various aspects of the path integral formulation on a firmer mathematical
footing \cite{AS92,AS94,BarNatan91,Kontsevich94,Rozansky97}.

Another strategy is to use the physics as a source of conjectures and of
geometric objects we should expect to see revealed in the combinatorial
structure if we look hard enough.  This is the strategy of the current paper.

The  geometry and the algebra part company in the first step of the
construction of the invariants, in which the geometric construction
begins with a 
compact semisimple Lie group, while the algebraic construction begins with a
semisimple Lie algebra.  These are almost, but not quite, in one-to-one
correspondence.  There are typically several Lie groups  with the
same Lie algebra, which differ only in their fundamental group.   The
invariants constructed from the Lie algebra correspond to the geometric
construction with the simply-connected Lie group.
Dijkgraaf and Witten \cite{DW90} address the issue of the existence of the
geometrically defined invariant for nonsimply-connected groups.

Let $G$ be a connected, simply-connected compact simple Lie
group with Lie algebra $\g,$ let $Z$ be a  subgroup of its center $Z(G),$
and let
$G_Z=G/Z$ be the quotient.  Recall $Z(G)=\Z_{l+1}$ if $G$ is of type
$A_l;$  $Z(G)=\Z_2$ if $G$ is of type $B_l,$ $C_l,$ or $E_7;$
$Z(G)=\Z_3$ if $G$ is of type $E_6;$ $Z(G)=\Z_4$ if $G$ is of type
$D_{2n+1};$ $Z(G)=1$ if $G$ is of type $E_8,$ $F_4,$ or $G_2;$ and
$Z(G)= \Z_2 \times \Z_2$ if $G$ is of type $D_{2n}.$      $G_Z$ like 
$G$ is  compact, simple
and connected with Lie algebra $\g,$ but its fundamental group  is
isomorphic to $Z.$  All connected 
compact simple Lie groups with Lie algebra $\g$ arise in this fashion.

Dijkgraaf and Witten argue that to construct a Chern-Simons theory
based on the group $G_Z$ requires only a choice of an element of
$H^4(BG_Z,\Z),$ that is  of the fourth cohomology of the classifying
space of $G_Z$ with integer coefficients.   Of course the projection map
from $G$ to $G_Z$ induces a homomorphism from $H^4(BG_Z,\Z) $ to
$H^4(BG,\Z).$  In every case but one (see below) these cohomology groups
are isomorphic to the integers, and thus the above homomorphism can be viewed
as multiplication by an integer $N.$  Dijkgraaf and Witten show that  $N$
 is the least $N$ such that $N(\lambda,\lambda)/2$ is an integer for
each fundamental weight $\lambda$ corresponding to an element of $Z,$ where
$(\,\cdot\,,\,\cdot\,)$ is the inner product in the weight space. We will find
it most convenient to index everything by the
\emph{level} $k \in H^4(BG,\Z),$ viewed as an integer, and
thus Dijkgraaf and Witten's work  predicts a
Chern-Simons theory associated to $G_Z$ exactly when $k$ is a multiple
of $N$ defined above.

 The one exception to the above observation about $H^4(BG_Z,\Z)$ is  when
$G$ is the simply-connected group associated to the Dynkin
diagram $D_{2n}$ and $Z$ is all of the center $\Z_2 \times \Z_2.$
We will follow Dijkgraaf and Witten in not considering this case,
although it should be extremely interesting (see, e.g., Felder,
Gawedski and Kupianen
\cite{FGK88})  and warrants further
study.  

The problem is then to construct the invariants from the quantum group
perspective.  Let us first briefly review how the Lie algebra appears in the
construction due to Reshetikhin and Turaev \cite{RT90,RT91}.

Beginning with a Lie algebra $\g,$ one deforms the universal enveloping
algebra $U(\g)$ by a deformation depending on a complex parameter $q$ to get an
algebra $U_q(\g)$  (called a quantum group) which satisfies the axioms of a
ribbon Hopf algebra.   In practice these axioms mean its representation theory can be
used to construct a system of link invariants.

For generic $q$ the algebra $U_q(\g)$ is semisimple, but when $q$ is a root of unity it
becomes nonsemisimple and quite subtle.  Unfortunately, the three-manifold
invariants arise only at roots of unity: In fact the Witten invariants at
level
$k$ discussed above correspond to $U_q(\g)$ when $q= \exp(2 \pi i/(k+h)),$
with $h$  the dual Coxeter number of $\g.$

Fortunately, the invariants do not depend on the quantum group directly, but
only on a piece of its representation theory.   In fact each of the
representations we consider corresponds naturally to a representation of the
original Lie algebra, and all the information we will need will be computed
from the classical representation using classical data involving weights and 
root spaces.   Any subset of this collection of representations satisfying
certain properties  (that it forms a \emph{modular category}) can be used
following the procedure of Reshetikhin and Turaev \cite{RT91,Turaev94}
to construct a  three-manifold invariant satisfying certain
cut-and-paste axioms expected of  topological quantum field theories 
(invariants  satisfying  these cut-and-paste axioms are called TQFTs in the
literature). 

A strategy thus naturally presents itself.  Nonsimply-connected Lie groups can
also be approached in terms of subsets of the set of representations of a Lie
algebra.  The finite-dimensional irreducible representations of $\g$
can be indexed by a cone inside the weight lattice called the Weyl chamber.   The
simply-connected group $G$ associated to $\g$ acts irreducibly on all these
representations, and in particular each element of the center $z \in  Z(G)$
acts on the representation indexed by a weight $\lambda$ as the identity times
the complex number $\chi_z(\lambda).$  In fact $\chi_z$ is a homomorphism from
the weight lattice to the unit circle for each $z,$ $\chi$ is a
homomorphism from $Z(G)$ into the dual group of the weight lattice, and for
each subgroup $Z$ the group $G_Z$ acts on exactly those representations whose weights
lie in the sublattice annihilated by $Z$ under $\chi.$  Thus our quantum
surrogate for $G_Z$ should be $U_q(\g)$ together with those representations
which lie in the sublattice annihilated by $Z.$  We have only to confirm that
this set of representations gives a modular category exactly when $k$ is a
multiple of $N$ as above.

	This strategy has been pursued in the special case of the
group $\rm{SO}(3),$ 
which is $\rm{SU}_2/(\Z_2),$ by Frohman and Kania-Bartoszy\'nskia
\cite{FK96} building on work of Kirby and Melvin \cite{KM91}. 
Unfortunately, while in this case Dijkgraaf and Witten predict a TQFT 
 and three-manifold invariant when $k$ is a multiple of $4,$ and in
fact suggest there should be some sort of spin TQFT and invariant  for
$k$ an odd multiple of $2,$  what
Frohman and Kania-Bartoszy\'nskia actually get is a modular category exactly when $k$
is odd!  We will see  that this holds for a general Lie
group: The appropriate set of representations forms a modular category 
only when
$k$ is relatively prime to Dijkgraaf and Witten's $N.$  It is difficult to
imagine a worse failure of the geometric predictions.

In fact there are two subtleties which  bring the
algebraic invariants  into line with the geometric predictions, though
there is  some additional structure in the situation for quantum groups which
is not readily apparent in the geometric point of view.  The first subtlety is
that in the cases where modularity is predicted but fails the failure is 
 because of a trivial sort of redundancy in the category which
can be readily quotiented out.
The result of the quotient is a modular category, and hence  a TQFT
and three-manifold invariant. 

 The second subtlety is that
many of the invariants coming from this quantum group construction can be
factored as a product of invariants (in fact the factoring works at the level
of TQFTs), one of which is a very simple invariant of the first homology
studied by Murakami, Ohtsuki and Okada \cite{MOO92} and the other of which
is the invariant associated to a smaller set of representations.    In the end
there is one prime invariant  (in the sense that it admits no further
decomposition as factors) for each TQFT conjectured by Dijkgraaf and Witten,
and all other invariants that we construct are formed out of these.   It is in
this sense that the conjectures of Dijkgraaf and Witten are confirmed.  
Interestingly, in many cases the invariants of Dijkgraaf and Witten  and the
original quantum group invariants constructed by Reshetikhin and Turaev are
not the prime version, but the prime invariant times one of the homology
invariants.   Since the homology invariants of Murakami et al (and hence the
invariants of Dijkgraaf and Witten)  are often zero this means that on closed
manifolds the prime invariants contain more topological information than those
that seem to arise naturally in the geometric interpretation.  On manifolds
with boundary the prime theory can be recovered from the composite theory.  

The first section of this paper  reviews the work of Brugui\`eres and
M\"uger on constructing quotients of a ribbon category which are
modular (this is quite distinct from the restriction to the Weyl alcove via
the truncated tensor product of Reshetikhin and Turaev \cite{RT91} and
Andersen and Paradowski \cite{AP95}, which we take as our starting point). 
In general there is a certain subcategory whose objects, called degenerate, are
obstructions to modularity, and which must form the kernel of the quotient.  In
the case at hand, where the ribbon categories possess a
$*$-structure, the quotient is possible as long as all the degenerate objects
are even, which is to say a change in the framing of components labeled by
these objects does not change the invariant.  For us, the degenerate objects
will always be invertible under the  tensor product, and in fact form a cyclic
group, which makes the modular category, invariant and TQFT easy to describe
concretely in terms of the original category.  The description of the
TQFT requires a much deeper immersion into the work of M\"uger and
Brugui\`eres than the rest of the paper, and invokes in detail the
theory of TQFTs, and thus is far less self-contained than the rest of
the article.  For this reason, despite the centrality of the TQFT to
the study of these invariants, the description is relegated to an appendix.

The second section deals with the quantum groups at roots of unity and their
representation theory. It considers the Weyl alcove, the subset
of weights corresponding to the representations of the quantum group that 
we are concerned with, and in particular the isometries of the Weyl alcove
(the weight lattice sits naturally in a Euclidean space).   Of central importance are
the invertible elements, which form the orbit of the trivial weight under
 isometries of the Weyl alcove.  We prove that  all degenerate
objects in these categories are
 invertible, and we identify these degenerate objects and determine
when they are even.   Most of the work of the section involves a
careful understanding of the so-called truncated tensor product of
representations, and relies on a crucial formula of Andersen and Paradowski,
 generalizing a classical formula for the ordinary tensor product of
representations to the quantum case.   
Isometries  and invertible elements of the Weyl alcove were
used for very similar ends in a paper by Felder et al. 
\cite{FGK88} to address Wess-Zumino-Witten theory for nonsimply-connected Lie
groups.  The section gives  a complete analysis of when a TQFT can be
constructed from the representations of any quantum group associated
to any classical simple, connected, compact Lie group at any level.
In particular,  this section offers  a prof purely in the language of
quantum groups
  that the full set of representations forms
a modular category and hence a TQFT.  The proofs in the literature
\cite{Kirillov96,Turaev94} 
all rely essentially on a result of Kac and Peterson from the theory of affine
Lie algebras \cite{KP84}.  Brugui\`eres has also constructed
modular categories  associated to the group $\mathrm{PGL}_n$ using his results.

Section Three identifies under what circumstances the  categories and TQFTs
constructed factor into simpler theories, and identifies the factors.   The
results of this section were first suggested to the author by E. Witten in
private conversation.   The decomposition exhibited here was observed 
previously by Kirby and Melvin \cite{KM91} on the level of invariants for the
group $\mathrm{SU}(2),$ where at odd levels the theory is a product of the prime
theory  (which they call the $\mathrm{SO}(3)$ theory) and Murakami et al's invariant $Z_2.$  
Again the more technical work with TQFTs appears in the appendix.

\subsection*{Acknowledgements}
I would like to thank Scott Axelrod, John Baez, John Barrett, Dana
Fine, Victor Kac, Gregor Masbaum, Isadore Singer, David Vogan, Eric Weinstein, 
and Edward Witten
for helpful conversations and suggestions.  I would like to thank
Michael M\"uger for alerting me to his work and that of Brugui\`eres,
which greatly simplified and extended this paper.

\section{Semisimple Ribbon $*$-Categories and Quotients}

We will follow the notation of Kirillov  \cite{Kirillov96} in this
section, and quote the basic results  from that paper, but other good
sources for the theory of ribbon categories include books by Kassel and Turaev
  \cite{Kassel95,Turaev94}.  This section relies heavily on the work
of M\"uger \cite{Muger99}.   Similar results were obtained
independently by  Brugui\`eres \cite{Bruguieres99}, but we will follow
M\"uger,  whose language is more in line with the rest of this paper.

\subsection{Semisimple ribbon $*$-categories}

A \emph{rigid   monoidal category} is a category
$\Cat$ 
together with a bifunctor $\trunc\colon \Cat \times \Cat \to \Cat$ which is
an associative multiplication with an identity object and morphism
(we will assume our category has been `strictified'  as in MacLane
\cite{MacLane71}) and a notion of duals compatible with $\trunc$. The example
on which to base one's intuition is the category of finite-dimensional vector
spaces and linear maps, where 
$\trunc$ is ordinary tensor product and the duals are vector space duals.

 The
category is called a semisimple $*$-category if the hom sets are vector spaces
over $\C,$
with  composition and  $\trunc$ acting as bilinear operations,  if
there are direct sums and kernels  (i.e., idempotents in the $\hom$
space of an object always factor through a subobject), and there is an antilinear
involution $*$ between $\hom(\lambda,\gamma)$ and $\hom(\gamma,\lambda)$ which is
contravariant $(f\circ g)^*=g^* \circ f^*$), monoidal ($(f\trunc
g)^*=f^* \trunc g^*$),  positive ($f^* \circ f=0 \implies f=0$), and
consistent with the duality in the appropriate sense.
  Here the example to keep in mind is the
category of representations of a $C^*$-algebra whose morphisms are linear maps between
representations which commute with the algebra's action, and whose simple objects
are the irreducible representations.  The $\trunc$ and duality structure occur
naturally if the $C^*$-algebra is a Hopf algebra.

\begin{figure}[hbt]
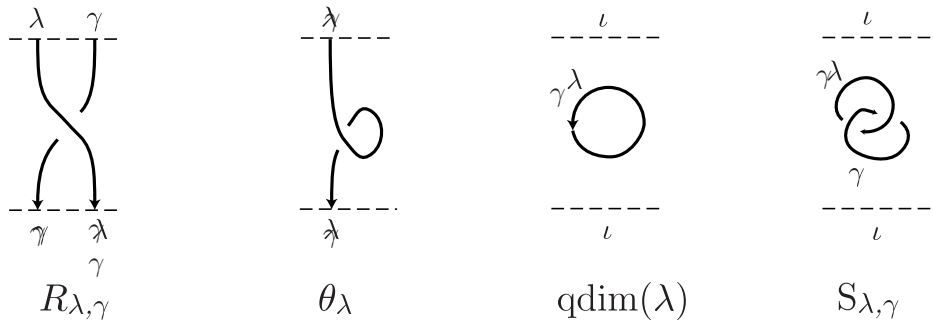
 
$$\pic{samples}{120}{-40}{0}{0}$$
\caption{Certain important knots and tangles and their image under the link
invariant} \label{fg:knots}
\end{figure}

 Finally, a semisimple  rigid monoidal
$*$-category $\Cat$ is a \emph{semisimple ribbon $*$-category} if every pair of
objects $\lambda$ and
$\gamma$ admits a unitary isomorphism $R_{\lambda,\gamma}\colon\lambda
\trunc \gamma \to \gamma \trunc \lambda$  and
every object $\lambda$ admits a unitary isomorphism
$\theta_\lambda\colon\lambda \to \lambda$ satisfying 
certain relations found in Kassel \cite{Kassel95}. In the presence
of the $*$-structure the $\theta$-morphism can be constructed from the
$R$-morphisms.   The relations are of
course designed to guarantee that there is a functor whose range is $\Cat$ and
whose domain is  the category of framed tangles with components labeled by
objects of
$\Cat$  \cite{Kirillov96,Sawin96a},  with a simple crossing  corresponding to
the
$R$ morphism and a full twist  to the $\theta$ morphism (see Figure
\ref{fg:knots}).   Note that the `quantum
dimension' $\qdim(\lambda)$ of each simple object $\lambda,$ which is defined in
\cite{Kirillov96} and corresponds to the invariant of the zero-framed
unknot labeled by $V$  (See Figure \ref{fg:knots}), is always a
positive real  (in fact $\geq 1$).  

Consider a semisimple ribbon $*$-category $\Cat.$  Let $\Gamma,$ the
\emph{label set} of $\Cat,$ be the set
of all isomorphism classes of simple objects in $\Cat.$  We will call
the  identity object for $\trunc$  $\iota,$ and will use the same name to
refer to its isomorphism class in  in $\Gamma.$  The dual
gives an involution on $\Gamma$ 
 which we  call $\dagger$.   Now for each $\gamma \in \Gamma$ the one-dimensional
$\hom(\gamma,\gamma)$ can be canonically identified with  $\C$ as a $C^*$- algebra,
so the morphism $\theta_\gamma$ corresponds to 
some complex number $C_\gamma$ of modulus $1$ such that
$\theta_\gamma$ is $C_\gamma$ times the identity.  Also if $\lambda,\gamma \in \Gamma,$ then
$\lambda \trunc \gamma$ is isomorphic to a sum $\bigoplus_{\eta \in \Gamma}
N_{\lambda,\gamma}^\eta \eta,$ where the nonnegative integers
$N_{\lambda,\gamma}^\eta$ represent multiplicities.   In
this and the sequel, we freely confuse simple objects with their
isomorphism classes, trusting the sophistication of the reader to
unravel the subtleties.   

Some facts about these numbers and their relation to the invariant
will be useful. We have $\qdim(\lambda^\dagger)=\qdim(\lambda),$
$\qdim(\lambda \trunc \gamma)= \qdim(\lambda)  \cdot \qdim(\gamma),$
$N_{\lambda, \gamma}^\eta=N_{\gamma, \lambda}^\eta=N_{\gamma,
\eta^*}^{\lambda^\dagger} = N_{\lambda^\dagger,
\gamma^\dagger}^{\eta^\dagger},$ $N_{\lambda,
\gamma}^\iota= \delta_{\lambda, \gamma^\dagger},$ and
$C_\lambda=C_{\lambda^\dagger},$
$C_\iota=1.$  

A ribbon category yields an invariant of labeled framed ribbon graphs.
 More specifically, consider an oriented graph embedded smoothly in
$S^3$ with a well-defined normal bundle  (that is, the edges incident
to a vertex are all tangent to a single plane at that point),
equipped with a nonzero section of the normal bundle and a choice of
edge at each vertex.  Label each
edge by an object of the category.  Notice the framing gives a cyclic
ordering on the edges incident to a vertex, and starting at the chosen
edge makes this a total ordering.   Thus to each vertex we can
associate the object $\overline{\lambda}_1 \trunc \overline{\lambda}_2 \trunc \cdots
\trunc \overline{\lambda}_n$ where $\lambda_i$ is the label of the $i$th incident
edge and $\overline{\lambda}_i$ is $\lambda_i$ or $\lambda_i^\dagger$ according to whether the
edge is oriented towards or away from the vertex.  Label each vertex
by an element of $\hom(\gamma,\iota),$ where $\gamma$ is this associated
object.   The category associates to this framed, labeled graph a number which
is invariant under ambient isotopy of the framed graph.  If a connected
component of the graph contains only bivalent vertices and thus is
homeomorphic to a circle (we call such components \emph{link components}) and
the edge label is simple,  we
can ignore the vertices and view it as a circle labeled by an object
of the category (up to an overall scale factor).  Thus we get
in particular a framed labeled link invariant.  The following
properties of the invariant will be important to us:

\begin{enumerate}

\item the invariant of a graph with an edge labeled by $\lambda \oplus
\gamma$ is the sum of the invariants of the same graph with that edge
labeled by $\lambda$ and $\gamma$ respectively, the labels on the
adjacent vertices being projected appropriately,
\item if the label of an edge is replaced by an isomorphic object and
the labels of the adjacent vertices are composed with the isomorphism
in the obvious way, the invariant is unchanged.   In particular, link
components can be unambiguously labeled by elements of $\Gamma,$ rather
than objects,
\item  the invariant of a graph with an edge labeled by $\iota$ is
the same as the invariant of the graph with that edge  deleted,
\item  the
invariant of a graph with an edge labeled by $\lambda$ is the
invariant of the graph with the orientation of that edge reversed
and the label replaced by $\lambda^\dagger,$ the labels of the adjacent
vertices remaining the same,
\item  the invariant of a graph with
a link  component labeled by $\lambda \trunc \gamma$ is the invariant of
the graph with that component replaced by two parallel components
(according to the framing) labeled by $\lambda$ and $\gamma$
respectively,
\item The invariant of the connected sum of two graphs  
along  edges labeled by a simple object $\lambda$ 
 is the product
of the invariants of the two graphs divided by $\qdim(\lambda)$, where
the connect sum consists of cutting the chosen 
edges and reconnecting them as in the definition of connected sum for knots, and
\item if a sphere intersects a ribbon graph only transversely at two
edges labeled by simple objects and oriented oppositely, the two
objects must be isomorphic or the invariant is zero.
\end{enumerate} 

We say a subset  $\Gamma'$ of $\Gamma$ is \emph{closed} if it is
closed under the duality involution and if whenever $\lambda, \gamma
\in \Gamma'$ and $N_{\lambda,\gamma}^\eta \neq 0,$ we have $\eta \in
\Gamma'$  (i.e., the product  of elements of  $\Gamma'$
is a sum of elements of $\Gamma'$).

\begin{prop}
If $\Gamma'$ is a closed subset of $\Gamma,$ the full subcategory of
$\Cat$ whose objects are sums of objects in the isomorphism classes in
$\Gamma'$ is again a semisimple  ribbon category.
\end{prop}

\begin{proof}
Immediate from the definition.
\end{proof}

\subsection{Degenerate objects}

Suppose $\Cat$ is a semisimple ribbon $*$-category with label set $\Gamma.$
For each $\lambda, \gamma \in \Gamma$ define $S_{\lambda,\gamma}$ to
be the value of the invariant of the zero-framed Hopf link with
components labeled by $\lambda$ and $\gamma$ respectively (see Figure
\ref{fg:knots}).  Thus
\begin{equation} \label{eq:S}S_{\lambda,\gamma}= \sum_\eta 
N_{\lambda,\gamma}^\eta    \qdim(\eta)
C_\eta C_\lambda^{-1} C_\gamma^{-1}.\end{equation}
By Properties 1-7 above 
$S_{\lambda,\gamma}=S_{\gamma,\lambda}=S_{\lambda^\dagger,\gamma^\dagger},$ 
$S_{\lambda,\iota}=\qdim(\lambda),$  $S_{\mu\trunc
\lambda,\gamma}=\sum_{\eta}N_{\mu,\lambda}^\eta S_{\eta,\gamma}$ and
$S_{\mu\trunc
\lambda,\gamma}=S_{\mu,\gamma}S_{\lambda,\gamma}/\qdim(\gamma).$    The
matrix of numbers
$S_{\lambda,\gamma}$ is called the $S$-matrix. Recall that a \emph{modular 
category} is a semisimple ribbon category  with a finite label set $\Gamma$
such that the $S$-matrix is invertible (see \cite{Kirillov96}).

M\"uger studies objects $\lambda$ such that
$R_{\lambda,\gamma}=R^{-1}_{\gamma,\lambda}$ for all $\gamma$ in the
category, and calls such objects \emph{degenerate}.   He proves

\begin{thm}\label{th:modular}
A semisimple ribbon $*$-category is modular if and only if all its
degenerate simple objects are isomorphic to the trivial object.
\end{thm}

Suppose one is computing the link invariant from a projection of a link with one
component labeled by a degenerate object.   By the definition, the
link invariant is unchanged by switching any of the crossings in which
that component participates from an overcrossing to an undercrossing
or vice-versa.   It is a standard observation of knot theory that in
that case the component can be unlinked from the other components and
unknotted except for the framing.  Thus the invariant of the original
link is the same as the invariant of the link with that component
deleted times $\qdim(\lambda) C_\lambda^m,$ where $m$ is the framing
of that component.  Moreover, from the fact that
$R_{\lambda,\lambda}^2=1$ we conclude that $C_\lambda=\pm 1.$  If
$C_\lambda=1$ we call $\lambda$ \emph{even} and we see that $\lambda$
provides no link information, and as far as the link invariant is
concerned behaves as if it were a multiple of the trivial object.   

One might reasonably suspect that one can `quotient out' by the even
degenerate objects to get a smaller ribbon category that in some sense
is a minimal ribbon category associated to the link invariant, and
further, that by doing so  one would make
the resulting quotient modular.  The first suspicion is correct, and
the second is correct as long as all degenerate objects are even.
M\"uger's proof of these results uses the fact that the degenerate
objects form a symmetric $*$-subcategory and the result of Doplicher
and Roberts that such a category is always the representation category
of a compact group.  Brugui\`eres uses the related Tannaka-Krein
theorem to prove a similar result which does
not assume the $*$-structure but replaces it with  assumptions that
amount to saying degenerate objects can be represented as vector
spaces.  The universal characterization of the quotient below appears in
Brugui\`eres' work and not in M\"uger's, but it is an easy consequence
of results in the latter.

\begin{thm} \label{th:quotient}
If $\Cat$ is a semisimple ribbon $*$-category such that all of its
degenerate objects are even, then there exists a modular category
$\Cat'$ and a ribbon $*$-functor from $\Cat$ to $\Cat',$ with the
property that every ribbon $*$-functor from $\Cat$ to a modular category factors
through this functor.
 In particular the link invariant of a link with
components labeled by objects of $\Cat$ is the same as that for the
link labeled by their image objects in $\Cat'.$
\end{thm}

\subsection{Invertible objects}

In general, the quotient category is quite complicated, and its
relation to the original category murky.   But in the special case
when all of the degenerate objects are invertible, which is the case for all
ribbon categories arising from quantum groups, the relationship can be
described much more explicitly.

An element $\lambda \in \Gamma$ is called  \emph{invertible} if $\lambda
\trunc \lambda^\dagger=\iota.$  
The set of invertible objects form a group under
tensor product, and if $\lambda$ is invertible and $\gamma$ is simple then
$\lambda \trunc \gamma$ is also simple, because its  product
with $\lambda^*$ is simple.  Thus each invertible object $\lambda$ corresponds to a map
$\phi_\lambda$ on objects defined by $\phi_\lambda(\gamma)= \lambda
\trunc \gamma$ 
which descends to a bijection on $\Gamma.$   The
map $\phi_\lambda$ satisfies the relation 
\begin{equation} \label{eq:phi}\phi_\lambda(\gamma \trunc
\gamma')=\phi_\lambda(\gamma) \trunc \gamma'.
\end{equation}
The set of maps $\phi_\lambda$ on $\Gamma$ forms a group under
composition  isomorphic to the group of invertible elements. 
We should caution that our confusion of objects and isomorphism
classes here generates a minor subtlety:  Isomorphic invertible elements
 generate distinct  maps of objects,
but all these maps descend to the same map on $\Gamma,$ which we
associate to the isomorphism class of the original objects.

Suppose $\Cat$ is a semisimple ribbon $*$-category all of whose
degenerate objects are both even and invertible.  The set of isomorphism
classes of such
objects forms an abelian group $Z$ under $\trunc.$ Suppose further
that this group $Z$ is cyclic,  so that in particular its second group
cohomology is trivial.   $Z$ acts on
$\Gamma$ by $\trunc,$ and associated to each element of $\Gamma$ is
its orbit and the stabilizer subgroup of $Z.$  M\"uger shows that in
the functor of Theorem \ref{th:quotient}  all the simple objects of an orbit
get sent to the same object in $\Cat',$ and that this object is the
direct sum of  a set of simple objects in $\Cat'$ in one-to-one
correspondence to the stabilizer  (in general there is a multiplicity
depending on a $2$-cocycle, but because of the cyclicity of the group,
this is one).  

In any ribbon category with finitely many isomorphism classes of simple objects we can
compute the  very important link invariant 
$$I(L)=\sum_{\gamma_1,\ldots, \gamma_n \in \Gamma}
\prod_{i=1}^n
\qdim(\gamma_i) F_{\gamma_1,\ldots, \gamma_n}(L)
$$
where $L$ is a link with $n$ components, the sum is over all ways of choosing $n$
isomorphism classes of simple objects out of $\Gamma,$ and $ F_{\gamma_1,\ldots,
\gamma_n}(L)$ is the invariant of $L$ with the $n$ components labeled by
$\gamma_1,\ldots,\gamma_n$ respectively  (or more properly, representative objects of each
of those classes).  If we extend the invariant to allow formal linear combinations of
objects as labels  (the invariant being linear on each label) this is the invariant of $L$
with each component labeled by $\omega=\sum_{\gamma \in \Gamma}\qdim(\gamma)\gamma.$

The importance of $I$ is that in the case of a modular category, 
 if $L$ is a surgery presentation of a $2$-framed three-manifold $M$ with $n$ components and
a linking matrix with signature $\sigma$ then
\begin{equation}\label{eq:invt}I(L)/I(H)^{n/2}
\end{equation}
depends only on $M$ and not on $L,$ where $H$ is the Hopf link.  Furthermore 
\begin{equation}\label{eq:finvt}(I(N)/I(P))^{\sigma/2} I(L)/I(H)^{n/2}
\end{equation}
depends only on the underlying three-manifold and not on the
$2$-framing, where $P$ is the $+1$ framed unknot and $N$ is the $-1$
framed unknot.  This can be written in the more congenial form
$$\left(\frac{I(N)}{|I(N)|}\right)^\sigma I(L)/|I(N)|^n$$
using the fact that $I(H)=I(P)I(N)$ and in a $*$-category $I(P)=\overline{I(N)}.$

\begin{prop}
Let $\Cat$ be a ribbon $*$-category such that all of its degenerate objects are even and
invertible, and such that the associated abelian group $Z$ is cyclic.  Then
$I(L)=|Z|^nI'(L),$ where
$I'$ is the corresponding invariant for the quotient category $\Cat',$ and therefore
Equations  (\ref{eq:invt}) and (\ref{eq:finvt}) give $2$-framed and ordinary three-manifold
invariants.
\end{prop}

\begin{proof}
Notice all simple objects in one orbit give the same link invariant,
so we could as well compute $I$ by taking the sum only over a
representative $\gamma$ of each orbit class, and replacing the factor
$\qdim(\gamma)$ with $|Z|/|S_\gamma|\qdim(\gamma),$ where $S_\gamma$
is the stabilizer of $\gamma$ and thus $|Z|/|S_\gamma|$
is the number of elements in the orbit of $\gamma.$  Because the
functor is a ribbon functor we can replace the invariant associated
with $\Cat$ with the invariant associated with $\Cat',$ if we replace
the label $\gamma$ on $L$ with its image under the functor.  This
image is the direct sum of $|S_\gamma|$ many simple objects, each with
quantum dimension $\qdim(\gamma)/|S_\gamma|.$  Decomposing the
invariant into a sum of invariants, each with the link labeled by a
single one of these simple objects, we obtain $|Z|^n I'(L).$
\end{proof}

\begin{rem}
    
   The result above should hold more generally.  In fact it is easy to 
   check that in any ribbon *-category all of whose degenerate objects 
   are even, the quantities (\ref{eq:invt}) and (\ref{eq:finvt}) are 
   respectively $2$-framed and ordinary three-manifold invariants.  It 
   is to be expected that they are equal to the corresponding 
   quantities in the quotient, but the proof may be more difficult.
\end{rem}

Of course we wish to describe the entire TQFT in
terms of the original category, just as we have here described the
three-manifold invariant.  In principle this is straightforward, but
in practice to do this concretely requires much more detail both about
M\"uger's construction and the construction of TQFTs from modular
categories than we will use in the rest of the article.   These
details appear in Appendix \ref{A:TQFT}, where a precise
description of the TQFT is given.

\begin{rem}
The term quotient which we use to denote the modularization is very appropriate
to the language of link invariants and TQFTs, but is not entirely accurate on
the level of the category.  In fact M\"uger proves that the functor from $\Cat$
to $\Cat'$ is an equivalence from $\Cat$ to $(\Cat')^{\hat{Z}},$ the subcategory
of
$\Cat'$ invariant under the action of the dual group to $Z.$
\end{rem}

\section{The Geometry of the Weyl Alcove}

This section relies heavily on the work of Andersen and Paradowski
\cite{AP95} and  of Kirillov \cite{Kirillov96}, 
which is summarized here.

\subsection{Quantum groups and the Weyl alcove}

Let $\g$ be a simple Lie algebra with Dynkin diagram different from
$D_{2n},$  and let $\uqg$ be its quantized
universal enveloping 
algebra.  This is defined exactly as in Kirillov with our $q$
equal to the square of his $q$
except we will normalize the inner product on
the Lie algebra so as to give long roots length $2$  (he normalizes so that
\emph{short} roots have length $2$).   This normalization has the properties
that it agrees with Kirby-Melvin and Reshetikhin-Turaev  
\cite{KM91,RT91} for $U_q(\mathrm{su}_2)$ and that it makes the 
quantum group $\uqg$ a modular Hopf algebra with the standard set of
representations exactly when $q$ is a root of unity.

We will need some notation from Lie algebra theory, most of which is taken from
Humphreys \cite{Humphreys72}, an excellent general reference on the subject.   Let
$r$ be the rank of $\g$ and let
$\{\alpha_i\}_{i \leq r}$ be the simple roots of $\g.$ The weight
lattice $\Lambda$ contains the sublattice $\Lambda_r$ spanned by the
roots, and we will be especially concerned with subgroups of the
fundamental group $\Lambda/\Lambda_r,$ because each such subgroup  corresponds  
to  a sublattice 
containing $\Lambda_r.$ The center $Z(G)$ of the simply-connected
group $G$ with Lie algebra $\g$ imbeds via the map $\chi$ defined in the
introduction into the dual group to $\Lambda,$ and in fact since each
element of the center acts trivially on the representations in the root
lattice, 
$\chi$ descends to an isomorphism from $Z(G)$ to the dual group of the
fundamental group $\Lambda/\Lambda_r,$ which isomorphism
we will also denote by $\chi.$

  The Weyl group is
denoted by $\W,$ and the set of weights in the fundamental Weyl chamber is
called
$\Lambda^+$  (we will loosely refer to this set itself as the Weyl chamber).  
Half the sum of the positive roots is called $\rho$  (Humphreys calls
this $\delta$), and  the unique long root in the Weyl chamber is called $\theta.$
This root corresponds to the adjoint representation of $\g.$  The dual Coxeter
number $h$ is defined to be $(\rho,\theta)+1,$ the value of the quadratic
Casimir on the adjoint representation.  The \emph{fundamental weights}
$\{\lambda_i\}_{i\leq r}$ are given by $(\lambda_i,\alpha_j)=
\delta_{i,j}(\alpha_i,\alpha_i)/2.$

Let $q=e^{2\pi i/(k+h)},$ for some natural number $k$. 
Kirillov shows that the category of representations of the quantum
group $\uqg$ corresponding to the weights in the \emph{Weyl alcove}
$\Lambda_0,$ i.e. those
$\lambda$ in the Weyl chamber such that $(\lambda,\theta)\leq k,$ form
a  semisimple ribbon category if the ordinary tensor product is
replaced by the 
\emph{truncated tensor product,} $\trunc,$ 
which returns  the  the ordinary tensor product 
 quotiented by the maximal tilting submodule \cite{AP95}.  
 Considered as a multiplication on the additive group of isomorphism
classes of representations spanned by those in the Weyl alcove  (with
direct sum as addition) it is  commutative, associative,
 distributive  and    determined by
\begin{equation} V_\lambda \trunc V_\gamma \iso
\bigoplus_{\eta \in \Lambda_0} N_{\lambda,
\gamma}^\eta
V_\eta,\end{equation}
where $N_{\lambda, \gamma}^\eta$ are nonnegative integers representing
multiplicities.  The principal result we use from Andersen and Paradowski
\cite{AP95} is their formula for these numbers, a variation on Racah's 
formula for  tensor product of classical representations
\begin{equation}  \label{eq:AP}
N_{\lambda,\gamma}^\eta = \sum_{\sigma \in
\W_0} (-1)^\sigma m_\gamma(\lambda-\sigma(\eta))
\end{equation}
where $m_\lambda(\mu)$ is the dimension of the $\mu$ weight space inside the
classical representation of highest weight $\lambda,$  $\W_0$ is the
quantum Weyl group, which is 
generated by  reflection about the hyperplanes
$\{x|(x+\rho,\alpha_i)=0\}$ for each simple root  $\alpha_i$ together with $\{x| 
(x,\theta)=k+1\},$  and $(-1)^\sigma$ is plus or minus one, according
to whether $\sigma$ is an even or odd product of these simple
reflections.  Notice that the weight $0,$ representing the trivial
representation, is the identity object for the truncated tensor product, and
thus is what in the last section was referred to by $\iota.$ Also
\begin{equation}\label{eq:theta}
C_\lambda=q^{(\lambda,\lambda+2\rho)/2}.\end{equation} 
and $\qdim(\lambda)>0.$

\subsection{The invertible elements of the Weyl alcove} 
\label{sec:invertible}

The geometric definition of the truncated tensor product  in
(\ref{eq:AP}) gives the invertible elements a special geometric role.  We will see
in the next subsection that degenerate simple objects are always invertible,
and thus that invertible elements will be of  central importance in the sequel.

\begin{thm}  \label{th:extreme} 
There is an injection $\ell$ from $Z(G)$ to the fundamental weights of the 
Weyl alcove such that $z$ acts on the classical representation 
$V_\gamma$ as $\exp(2 \pi i (\gamma,\ell(z))) \cdot \operatorname{id}_\gamma.$  The 
fundamental weights $\lambda_i$ in the image of $\ell$ are exactly those 
for which $(\lambda_i,\theta)=1$ and the associated root $\alpha_i$ is 
long, and are also exactly those for which there is a unique element $\tau_i$ 
of the classical Weyl group taking the standard base to the base 
$\{\alpha_j\}_{j \neq i} \cup \{-\theta\}.$  If we 
define $\phi_i(\gamma)=k\lambda_i +\tau_i(\gamma),$ then $\phi_i$ is an 
isometry of the Weyl alcove and of the simplex 
$\{\lambda:(\lambda,\alpha_j) \geq 0 \text{ and } (\lambda,\theta)\leq k\}$ 
and $\phi_i(\lambda \trunc \gamma)=\phi_i(\lambda) \trunc \gamma$ (i.e., 
$\phi_i$ satisfies Equation (\ref{eq:phi})).   If we use $k$ also to 
represent the map on the weight lattice which multiplies each weight 
by the 
number $k,$ then $k\ell$ is a homomorphism in the sense that 
$k\ell(zz')=k\ell(z) \trunc  k\ell(z').$  Weights in the range of $k\ell$ can be 
characterized as extreme points of the simplex  $\{\lambda:(\lambda,\alpha_j) 
\geq 0 \text{ and } (\lambda,\theta)\leq k\}$ such that a neighborhood of 
the weight $0$  in the simplex is isometric to a neighborhood 
of the extreme point in the simplex, the isometry being given 
by $\phi_i.$
\end{thm} 

\begin{figure}[hbt]
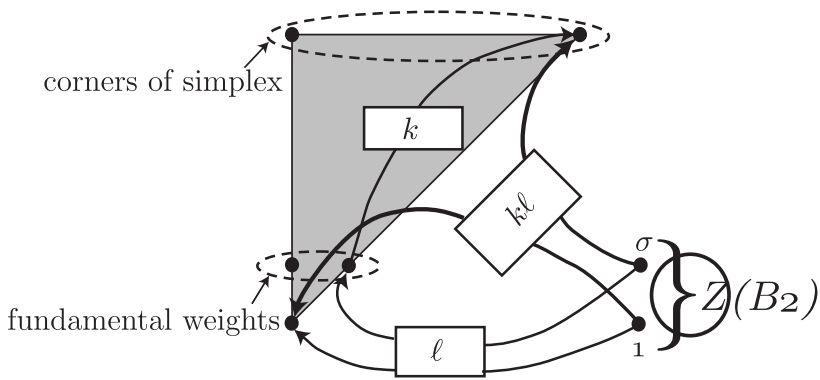
 
$$\pic{example}{150}{-60}{0}{0}$$
\caption{The maps $\ell,$ $k,$ and $k\ell$ for the case $B_2$} 
\label{fg:example}
\end{figure}

\begin{rem}  \mbox{}
\begin{itemize}
\item  This
$\ell$ is the isomorphism referred to in the introduction.
Specifically,  Dijkgraaf and Witten predict Chern-Simons theories
for $G_Z$ when  $Z$ is a subgroup of $Z(G)$ and $k$ is
such that $k(\ell(z),\ell(z))/2$ is an integer for each
$z \in Z.$
\item The homomorphism property of $k\ell$ shows its range consists of 
invertible objects, and of course $\phi_i$ is the map satisfying 
Equation   (\ref{eq:phi})  which  is associated to 
the invertible object $k\lambda_i.$    It 
is almost the case that these are all the invertible elements of the Weyl 
alcove.  In fact it is shown in \cite{Sawin00c} that the only case in which 
these are not all of the invertible elements is for the quantum group of 
type $E_8$ at level $k=2,$ when the fundamental weight $\lambda_1$ is 
invertible but not in the range of $k\ell.$  In fact, in this case the 
associated map $\phi$ is an isometry of the Weyl alcove, but 
\emph{not} of the simplex.
\item The local isometry condition says more intuitively that the range of 
$k\ell$ is the set of `sharp corners' of the simplex.  from the 
proposition 
it follows that the group of isometries of the simplex is the 
semidirect product of the  group of isometries of the Weyl chamber (which 
correspond naturally to automorphisms of the Coxeter graph) with the group 
of maps $\{\phi_i\}.$
\item Of course more generally there is a bijection from the full set of 
fundamental weights to the full set of extreme points  (or corners) of the 
simplex given by $\lambda_i \mapsto k\lambda_i/(\lambda_i,\theta).$  For 
appropriate $k$ these corners are weights and will figure prominently in 
the next section.
\item The maps $\ell,$ $k,$ and $k\ell,$ their ranges and domains and their 
relation to the geometry of the simplex are illustrated in the case of 
$B_2$ for some arbitrary $k$ in Figure \ref{fg:example}.   The map 
$\phi_i$ 
associated to the nontrivial element $\sigma$ of $Z(B_2)=\Z_2$ is of course 
the reflection about the diagonal.
\end{itemize}
\end{rem}

\begin{lem}  \label{lm:extreme} 
The set of fundamental weights $\lambda_i$ such that $\alpha_i$ is long 
and $(\theta,\alpha_i)=1$ is the same as the set such that there is a unique 
element $\tau_i$ of the classical Weyl group taking the standard base 
to the base $\{a_j\}_{j \neq i} \cup \{-\theta\}.$  
If $e$ is the homomorphism  from $\Lambda$ to $(\Lambda/\Lambda_r)^*,$ 
where $(\Lambda/\Lambda_r)^*$ is the dual 
group to the weight lattice modulo the root lattice, which sends each 
weight $\lambda$ to the homomorphism $e_\lambda$ given by 
$e_\lambda(\gamma)=\exp(2 \pi i (\lambda,\gamma)),$ then $e$ is a bijection when 
restricted to this set.
\end{lem} 

\begin{proof}
See the end of Section \ref{sec:invertible}.
\end{proof}

\begin{proof}[Proof of Theorem \ref{th:extreme}]

It is well known  that $Z(G)$ is isomorphic to the 
group $(\Lambda/\Lambda_r)^*$ by the map sending $z$ to the homomorphism 
$\chi_z$ on $\Lambda$ such that $z$ acts on the classical representation 
$V_\gamma$ by $\chi_z(\gamma) \operatorname{id}_\gamma$  (this is clearly a homomorphism 
descending to $(\Lambda/\Lambda_r)^*,$ and is injective by  the 
faithfulness of the left regular representation.  That the domain and range 
have the same dimension is shown in \cite{Humphreys72}).  Thus by Lemma 
\ref{lm:extreme} we can construct a bijection $\ell$ from $Z(G)$ to the 
set of fundamental weights $\lambda_i$ meeting the two characterizations of 
the proposition.   This gives the first two sentences of the proposition.  For the 
rest, we argue as follows.  We first show that the map $\phi_i$ defined in 
the proposition satisfies the conditions of the proposition.  We then use 
this to give the characterization of the range of $k\ell$ in the last 
sentence of the proposition.  Finally this will allow us to show that 
$k\ell$ is a homomorphism.

The map $\phi_i$ is certainly an isometry of the weight space taking
the weight $0$ (i.e. the object $\iota$) to $k\lambda_i.$   The image of the elements of the
Weyl alcove are those weights $\gamma$ such that $(\gamma-k\lambda_i,\alpha_j)
\geq 0$ for $j \neq i,$ $(\gamma-k\lambda_i,-\theta) \geq 0,$ and
$(\gamma-k\lambda_i, -\alpha_i)\leq k.$  Since $(k\lambda_i,\alpha_j)=0$ for
$j\neq i$ and $(k\lambda_i,\theta)=(k\lambda_i,\alpha_i)=k,$ such $\gamma$ are
exactly the elements of the Weyl alcove.  Thus $\phi_i$ is an isometry 
taking the Weyl alcove (and the simplex)
to itself.

To see that $\phi_i(\lambda \trunc
\gamma)=\phi_i(\lambda) \trunc \gamma,$ note that for any element
$\sigma$ of the quantum Weyl group $\sigma' \circ \phi_i= \phi_i \circ
\sigma,$ where $\sigma'$ is another element of the quantum Weyl group
with the same sign, because this is true for generating reflections.  Thus in
Formula (\ref{eq:AP})
\begin{multline*}N_{\lambda,\gamma}^\eta  = \sum_{\sigma \in
\W_0} (-1)^\sigma m_\gamma(\lambda-\sigma(\eta))=\sum_{\sigma \in
\W_0} (-1)^\sigma m_\gamma(\tau_i(\lambda-\sigma(\eta)))
\\=\sum_{\sigma \in
\W_0} (-1)^\sigma m_\gamma(\phi_i(\lambda)-\phi_i(\sigma(\eta)))=
\sum_{\sigma' \in
\W_0} (-1)^{\sigma'}
m_\gamma(\phi_i(\lambda)-\sigma'(\phi_i(\eta)))
\\=N_{\phi_i(\lambda),
\gamma}^{\phi_i(\eta)}.
\end{multline*} This confirms $\phi_i(\lambda \trunc 
\gamma)=\phi_{i}(\lambda) \trunc \gamma.$

Now $\phi_{i}$ certainly takes a neighborhood of the weight $0$ to a 
neighborhood of the weight $k\lambda_{i},$ and since it is an isometry 
of the simplex it also connects these neighborhoods intersected with 
the simplex.   Conversely, if $\lambda$ is an extreme point of the 
simplex and $\phi$ is an isometry of a neighborhood of the weight $0$ 
intersected with the simplex to a neighborhood of $\lambda$ 
intersected with the simplex, then $\phi$ extends to an isometry from 
an entire  neighborhood of the weight $0$ to a neighborhood of the 
weight $\lambda,$ which takes the hyperplanes 
$\{\gamma:(\alpha_{j},\gamma)=0\}$ to the hyperplanes 
$\{\gamma:(\alpha_{j},\gamma)=0\}$ for all $j \neq i$ and for some $i$ 
together with the hyperplane $\{\gamma:(\theta,\gamma)=(\theta,\lambda)\}.$  
Therefore $\gamma \mapsto \phi_{i}(\gamma) -\lambda$ is an isometry 
taking $0$ to $0$ and the hyperplanes 
$$\{ \{\gamma:(\alpha_{j},\gamma)=0\}: 0 \leq j \leq r\}$$
to 
$$\{\{\gamma:(\alpha_{j},\gamma)=0\}:j \neq i\} \cup 
\{\{\gamma:(\theta,\gamma)=0\}\}.$$
Such an isometry is necessarily a composition $ \tau \tau'$  with 
$\tau$ an    element 
of the classical Weyl group and $\tau'$ a linear isometry permuting the simple roots  
$\{\alpha_{j}\}$   (such 
isometries of the Weyl chamber correspond to automorphisms of the 
Coxeter graph).  The map $\tau$ corresponds to the base 
$\{\alpha_j\}_{j \neq i} \cup \{-\theta\},$    and
therefore by Lemma  \ref{lm:extreme}, the associated $\lambda_{i}$ is 
in the range of $\ell.$  Since the composition of isometries $\phi 
\circ (\tau')^{-1}\phi_{i}^{-1}$ of the alcove takes $\gamma$ to 
$\gamma -k\lambda_{i} + \lambda,$ we must have that 
$\lambda=k\lambda_{i}$ and we see that every corner of the simplex 
which is locally isometric to the corner $0$ is in the image of 
$k\ell.$  

To see that $k\ell$ is a homomorphism, let us first consider the case 
$k=1.$   Then if $\ell(z_{1})=\lambda_{i}$  we have $\ell(z_{1}) 
\trunc \ell(z_{2})= \phi_i(0) \trunc 
\ell(z_{2})=\phi_{i}(\ell(z_{2}))$ which is again an extreme point 
locally isometric to $0,$ and thus is $\ell(z_{3}) $ for some 
$z_{3}.$  Of course in general each $\eta$ such that 
$N_{\lambda,\gamma}^{\eta} \neq 0$ differs from $\lambda +\gamma$ by 
an element of the root lattice, and thus 
$e_{\eta}=e_{\lambda+\gamma}=e_{\lambda}\cdot e_{\gamma}.$  In 
particular $e_{\ell(z_{3})}=e_{\ell(z_{1}z_{2})}$ and 
$z_{3}=z_{1}z_{2}.$  Thus $\ell$ is an endomorphism.

Finally notice that the map multiplication by $k$ where it is  defined  forms a 
commuting  square with the maps $\phi_i$ acting 
respectively on  the Weyl 
alcove at level $1$ and  the Weyl alcove at level $k$  and thus the 
map $k$ takes on the one hand 
$\ell(z_{1}z_{2})$ to $k\ell(z_{1}z_{2})$ and on the other 
$\ell(z_{1}) \trunc \ell(z_{2}) = \phi_{i}(\ell(z_{2}))$ to 
$\phi_{i}(k\ell(z_{2}))=k\ell(z_{1}) \trunc k\ell(z_{2})$ and thus 
$k\ell$ is an endomorphism.
\end{proof}

\begin{figure}[hbt]
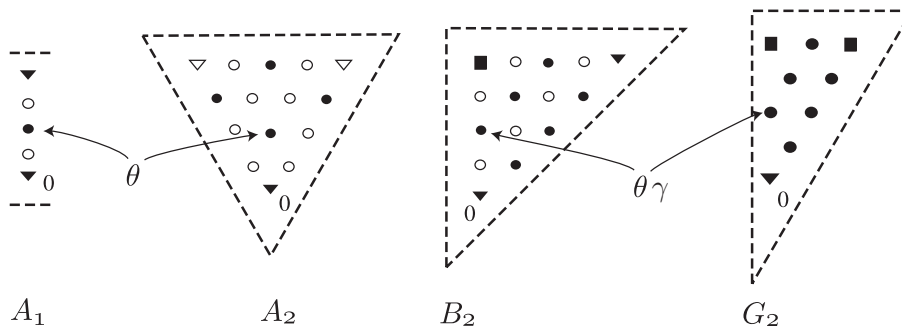


$$\pic{weyl}{120}{-55}{0}{0}$$

\caption{Weyl alcove and corners for low rank Lie algebras} \label{fg:units}
\end{figure}

Figure \ref{fg:units} illustrates the notions discussed above for rank one and
two Lie algebras.  The Weyl alcove for $k=4$ is shown.  Invertible weights are
marked with a triangle (including the weight $0,$ which is labeled), other
corners with a square.     The lines of reflection generating the quantum Weyl
group are  indicated by  dotted lines, and the root $\theta$ is marked.   The
elements of the root lattice are indicated by solid dots, squares and 
triangles and  weights not in the
root lattice by open figures.

\begin{proof}[Proof of Lemma \ref{lm:extreme}]

Suppose $\alpha_i$ is long, and $(\lambda_i,\theta)=1.$   Since
$(\lambda_i,\theta)$ is the coefficient of $\alpha_i$ in the expansion
of $\theta$ in terms of simple roots, we have 
$$\theta= \alpha_i + \sum_{j \neq i} k_j \alpha_j$$
for some positive integers $k_j$.   To conclude the existence of a unique
$\tau_i$ as in the statement of the lemma, it suffices to show that
$\{\alpha_j\}_{j\neq i} \cup 
\{-\theta\}$ is a base  Since the coefficient of $\alpha_i$
in the expansion of $\theta$ is $1,$ 
every positive root can be written in
terms of the original base with the coefficient of $\alpha_i$ being either $0$ or
$1.$   In the former case the root is already a nonnegative combination of
$\{\alpha_j\}_{j\neq i},$ in the second  $\alpha= \theta + (\alpha-
\theta)$ writes the positive root as a sum of two terms, each with all
nonpositive coefficients in the new base. 

That the second condition implies the first is a straightforward
inverting of the above argument.

For the second assertion, notice that the homomorphism is
well-defined  on $\Lambda/\Lambda_r,$ because $(\lambda_i,\alpha_j)=
\delta_{i,j},$ so $(\lambda_i,\alpha)$ is an integer for all roots
$\alpha.$   A count of such $\lambda_i$ from Table 2 of
\cite{Humphreys72}[Chapt. 12] (For $A_l$ all fundamental weights, for
$B_l$ $\lambda_1,$ for $C_l$ $\lambda_l,$ for $D_l$ $\lambda_1,$
$\lambda_{l-1}$ and $\lambda_l,$ for $E_6$ $\lambda_1,$ $\lambda_6,$
and for $E_7$ $\lambda_7$) shows that we need only check that no
$\lambda_i$ gets sent to the trivial homomorphism.  This is confirmed
by computing inner products $(\lambda_i,\lambda_j)$ using Table 1 of
\cite{Humphreys72}[Chapt. 13]. 
\end{proof}

\subsection{Degenerate objects are invertible} \label{sec:degenerate}

 Recall that $\theta$ represents the unique long root in the
Weyl chamber.   Let us use $\beta$ to represent the unique short root in the
Weyl chamber.   

\begin{prop} \label{pr:closed} 
Let $Z$ be a subgroup of $Z(G).$  Let $\Gamma_Z$ be the set of weights
in $\Lambda_0$ which are annihilated by $\chi_{z}=e_{\ell(z)}$ for 
all $z \in Z.$  Notice $\Gamma_Z$
is the sublattice of weights  corresponding to classical
representations of $G_Z,$ intersected with $\Lambda_0.$  Let
$\Delta_Z$ be $k\ell[Z].$  Then $\Gamma_Z$ and $\Delta_Z$ are
closed subsets of $\Lambda_0.$
\end{prop}

\begin{proof}
By Equation (\ref{eq:AP}), the weights occurring in the decomposition
of the truncated tensor product of two weights lie in the product of
their cosets in $\Lambda/\Lambda_r,$ and thus are annihilated by any
homomorphism which annihilates the factors (since the reflections that
generate the quantum Weyl group preserve these cosets).   Thus
$\Gamma_Z$ is closed under the truncated tensor product.  Of course it
is closed under the duality relation, since duality corresponds to
inverse in $\Lambda/\Lambda_r.$

In $\Delta_Z$ the truncated tensor product and dual  on simple objects
corresponds to product and inverse in $Z,$ so closure is immediate.
\end{proof}

\begin{rem}
The sets $\Gamma_{Z}$ we may view as classical closed sets.  The 
closed subsets of the Weyl chamber are the restriction of sublattices 
of $\Lambda$ containing the root lattice to the Weyl chamber, and the 
subsets $\Gamma_{Z}$ are just the restriction of these to the 
alcove.   the sets $\Delta_{Z},$ however, have no correspondence to 
anything classical.  These two classes of closed subset almost, but 
not quite, exhaust the list.   It is shown in \cite{Sawin00c} that 
when (and only when) $k=2$ there are certain additional closed subsets.
\end{rem}

To address the question of whether a subset of the form $\Gamma_Z$ or
$\Delta_Z$ yields a modular category via the quotient of Section 1, we need to identify the
degenerate objects of these sets.   Since we have an explicit
formula for $C_\lambda,$ this is largely a matter of using Equation
(\ref{eq:AP}) to give a careful description of the truncated tensor
product.

\begin{lem} \label{lm:exist} 
For any $\sigma$ in the classical Weyl group $\W,$ and any weights
$\gamma,$ $\lambda$ in the Weyl alcove, if
$\lambda+\sigma(\gamma)$ is in the Weyl alcove,
 then
$\lambda \trunc \gamma$ contains $\lambda+\sigma(\gamma)$ as a summand
with multiplicity one.
\end{lem}

\begin{proof}
See end of Section \ref{sec:degenerate}.
\end{proof}

\begin{lem} \label{lm:theta}
$\lambda \trunc \lambda^\dagger$ contains $\theta$ as a summand  if $k \geq
2$ and $\lambda$ is not a corner  (i.e. a multiple of a fundamental
weight such that $(\lambda,\theta)=k$).  In the nonsimply-laced case it 
contains $\beta$ as a summand  unless $(\lambda,\alpha_i)=0$ for every
short simple root $\alpha_i.$ 
\end{lem}

\begin{proof}
See end of Section \ref{sec:degenerate}.
\end{proof}

\begin{thm} \label{th:pseudo}
The set of degenerate objects of a closed subset of the form
$\Gamma_Z$ or $\Delta_Z$ is a set of the from $\Delta_{Z'}$ for some
subgroup $Z'$ of $Z(G).$  In particular, all degenerate objects are
invertible. 
\end{thm}

\begin{proof}
Of course since $\trunc$ on $\Delta_Z$ can be identified with group
multiplication on $Z,$ the only closed subsets of $\Delta_Z$ will be 
subgroups and hence of the
form $\Delta_{Z'}.$  So assume $\Gamma=\Gamma_Z,$ and that $\lambda
\in \Gamma$ is degenerate.  We will show that $\lambda$ is in the
range of $k\ell,$ which 
suffices for the theorem.

In that case if $k>1$ then $\theta, \beta \in \Gamma,$ so
$$C_\gamma C_\lambda^{-1} C_\theta^{-1}=1$$
for any $\gamma$ with $N_{\lambda,\theta}^\gamma \neq 0,$ and likewise
for $\beta.$

If $\lambda$ is not a corner of the Weyl alcove then
$N_{\lambda,\theta}^\lambda \neq 0$ by Lemma \ref{lm:theta} so
$$
C_\lambda C_\lambda^{-1} C_\theta^{-1}
=\exp\left(\frac{-\pi i[ (\theta,\theta)+2(\theta,\rho)]}{k+h}\right)
=\exp\left(\frac{-2 \pi i h}{k+h}\right)\neq 1
$$
so $\lambda$ is not degenerate.  Likewise in the nonsimply-laced
case if $(\lambda,\alpha_i) \neq 0$ for some short simple root
$\alpha_i$ then by Lemma \ref{lm:theta} $N_{\lambda,\beta}^\lambda
\neq 0$ so
$$
C_\lambda C_\lambda^{-1} C_\beta^{-1}
=\exp\left(\frac{-\pi i [(\beta,\beta)+2(\beta,\rho)]}{k+h}\right) \neq 1
$$
since $(\beta,\beta) + 2(\beta,\rho) < (\theta,\theta)+ 2(\theta,\rho)
= 2h,$ so $\lambda$ is not degenerate.

Now if $\lambda+\alpha$ is in the Weyl alcove for some long root
$\alpha$ then by Lemma \ref{lm:exist}
$N_{\lambda,\theta}^{\lambda+\alpha} \neq 0$ so
\begin{multline*}
C_{\lambda+\alpha} C_\lambda^{-1} C_\theta^{-1}
=\exp\left(\frac{\pi i[2 (\lambda,\alpha)+2(\alpha,\rho)-
2(\theta,\rho)]}{k+h}\right) \\ = \exp \left(\frac{2 \pi
i[(\lambda+\rho,\alpha)-(\theta,\rho)]}{k+h}   \right).
\end{multline*}
Also $|(\lambda+\rho,\alpha)| \leq k+h$ for all $\alpha,$ so
$(\lambda+\rho,\alpha)-(\theta,\rho)$ can only be a multiple of $k+h$
if 
\begin{enumerate}
\item $(\lambda+\rho,\alpha)=h-1$ or
\item $(\lambda+\rho,\alpha)=-(k+1).$
\end{enumerate}

If $\lambda$ is a corner and is orthogonal to all short simple roots,
then $\lambda$ is $k\lambda_i/n$ for $\alpha_i$ long, where
$(\lambda_i,\theta)=n.$  If $k>n$ then $\lambda-\alpha_i$ is in the
Weyl alcove,  $N_{\lambda,\theta}^{\lambda-\alpha_i} \neq 0,$ and 
$$(\lambda+\rho,-\alpha_i)=-k/n-1.$$
For $\lambda$ to be
degenerate requires this quantity to be equal to $-k-1$ (since it
is negative, it is not $h-1$), so that $n=1$ and we conclude $\lambda$
is in the range of $k\ell$ by Lemma \ref{lm:extreme}.

Thus the only possible degenerate objects which are not in the range 
of 
$k\ell$ are weights $\lambda_i$ dual to long roots for
$k=(\lambda_i,\theta).$ We argue first that for each such $\lambda_i$
there is a long positive root $\alpha$ with $(\lambda_i,\alpha)=0$ such that
either $\lambda_i + \alpha$ is in the Weyl alcove with
$k=(\lambda_i,\theta)$ or $\alpha-\theta$ is a long root and
$\lambda_i +\alpha-\theta$ is in the Weyl alcove, and second that this
contradicts degeneracy.

To see the existence of such an $\alpha,$ observe from the Dynkin
diagrams \cite{Humphreys72}[pg. 58] that for every fundamental weight
$\lambda_i$ dual to a long root, either $(\lambda_i,\theta)=1,$ 
$\lambda_i=\theta,$ or for one of the subdiagrams into which the removal of
$\lambda_i$ divides the diagram,  the weight $\lambda_j$ adjacent to
$\lambda_i$  is dual to a long root and satisfies
$(\lambda_j,\theta')=1$ for $\theta'$ the highest root associated to
this subdiagram.   In the first case $\lambda_i$ is in the image of  $k
\ell,$  in the second $\alpha=\iota$ will do, and in the third we
choose $\alpha=\theta'.$  Of course in the third case $(\alpha,\lambda_i)=0,$ 
$(\alpha,\lambda_k)\geq 0$ for $k \neq i,$ and $(\lambda_i +
\alpha,\alpha_i)=-1 + 1$ because the decomposition of $\alpha$ into
simple roots contains exactly one simple root  adjacent to $\alpha_i.$
So $\lambda_i+\alpha$ is in the Weyl chamber.  If
$(\alpha_i,\theta)=0,$  then $\lambda_i + \alpha$ is in fact in the
Weyl alcove.   If not then $(\alpha,\theta)=1.$  Except for the case
$A_l,$ where all corners are in the range of $k \ell$ and there is nothing
to prove, $\theta$ is of the form $\lambda_k$ for some $k,$ so
$(\alpha,\theta)=1$ indicates that the subdiagram  contained that $\lambda_k.$
Further inspection of the Dynkin diagrams indicates that the
$\lambda_i$ for which the only subdiagram meeting the desired
conditions  contains this $\lambda_k$ are $\lambda_2$ of $E_7$ and
$\lambda_1$ of $E_8.$  In the first case $(\alpha,\alpha_k)=1$ and
thus $\lambda_2+\alpha-\theta$ is in the Weyl alcove.  In the second
use $\alpha+\alpha_8$ is a positive root, $\alpha+\alpha_8-\theta$ is
a root, and $\lambda_i +\alpha+\alpha_8-\theta$ is in the Weyl alcove,
so $\alpha+\alpha_8$ meets the desired condition.

Thus we need only check that the existence of such $\alpha$
contradicts degeneracy.   In the first case
$$0<(\lambda+\rho,\alpha)=(\rho,\alpha)<h-1,$$
because of course $\alpha \neq \theta.$  In the second
$$0>(\lambda+\rho,\alpha-\theta)=-k +(\rho,\alpha)>-k-1.$$
Thus the
only degenerate objects for $\Gamma$ are those in $k 
\ell[Z(G)].$

Now if $k=1,$ it is not true that $\theta \in \Gamma.$  However, the
only elements of $\Lambda_0$ and hence of $\Gamma$ are elements of $k
\circ \ell[Z(G)]$ and weights $\lambda_i$ with $\alpha_i$ short.  Thus
the argument involving $C_\beta$ above suffices to show that only
those in $k \ell[Z(G)]$ can be  degenerate.
\end{proof}

\begin{cor}
Every closed subset of the form $\Gamma_Z$ or $\Delta_Z$ yields, via
the quotient of Section 1, a modular category whose associated TQFT
and invariant is as described in Section 1 and the Appendix.
\end{cor}

\begin{proof}
Since we have exempted $D_{2n},$ $Z$ is cyclic.
\end{proof}

\begin{proof}[Proof of Lemma \ref{lm:exist}]
 We note first that if $\lambda$ is in the Weyl alcove, and $\sigma$ is an
element of the quantum Weyl group taking a weight $\mu$ not in the Weyl alcove
into the Weyl alcove, then the distance between $\lambda$ and $\sigma(\mu)$ is
strictly less than the distance between $\lambda$ and $\mu.$  To see this, note
that if $\sigma(\mu)$ is in the Weyl alcove, than $\mu$ cannot lie on one of the
`walls of the Weyl alcove,' i.e. the hyperplanes reflection about which
generates the quantum Weyl group  (if it did, it and all its conjugates would
have nontrivial stabilizers, which is not true of any point in the Weyl
alcove).   In that case, one of the walls of the alcove lies between $\lambda$
and $\mu,$ and thus reflection about this wall brings $\mu$ strictly closer to
$\lambda.$  repeating this procedure brings a sequence of weights conjugate to
$\mu$ getting strictly closer to $\lambda.$  Since there are only finitely many
weights in the weight lattice a given distance from $\lambda,$ this process
must end after finite time.   It can only end by reaching a point which is
conjugate to $\mu$ and which lies in the alcove or on the walls.
Since the quantum Weyl group acts transitively on the 
Weyl alcove  \cite{KP84}, this gives the claim.

Now let $\lambda,$ $\gamma,$ and $\sigma$ be  as in the statement of
the lemma, so
any $\mu$ for which $m_\gamma(\lambda-\mu)$ is nonzero must be a distance at
most $||\gamma||$ from $\lambda,$ so if $\mu$ is not in the Weyl alcove but is
conjugate to $\mu'$ which is, then $\mu'$ is a distance less than $||\gamma||$
from $\lambda,$ and hence is not $\lambda+\sigma(\gamma).$   Thus the only
contribution to $N_{\lambda,\gamma}^{\lambda+\sigma(\gamma)}$ in Formula
(\ref{eq:AP}) comes from $\sigma=1,$ and would be $m_\gamma(\sigma(\gamma))=1.$
\end{proof}

\begin{proof}[Proof of Lemma \ref{lm:theta}]
  We will make the argument for $\theta,$ noting parenthetically how
it differs for $\beta$ when not simply-laced.  We will actually 
prove that $N_{\lambda,\theta}^\lambda$ is nonzero, which is
equivalent.  Recall that since $\iota$ is the weight $0$ $m_\theta(\iota)=r$
(respectively
$m_\beta(\iota)=r_0,$ the number of short simple roots).  Thus in the sum
(\ref{eq:AP}), there is a contribution of $r$ ($r_0$) from the identity element
of the quantum Weyl group and a 
contribution for each $\sigma$ such that $m_\theta(\lambda,\sigma(\lambda))
\neq 0.$  By the first paragraph in the proof of Lemma \ref{lm:exist}
above, we saw that 
$\sigma$ can be written as a product of reflections each taking $\lambda$
strictly farther away from itself.   Since $\sigma(\lambda)-\lambda$ is in the
root lattice, we conclude that after one such reflection its length is at least
that of a short root, after two its length is at least that of a long root, and
after three it must be longer than a long root.   Thus if $\sigma $ is such a
product of three or more reflections, $m_\theta(\lambda-\sigma(\lambda))=0.$  A
product of two reflections only increases $N_{\lambda,\theta}^\lambda,$ so it
suffices to consider the effect of a single reflection.   If $\sigma$ is
reflection about one of the  walls of the Weyl alcove and
$m_\theta(\lambda-\sigma(\lambda)) \neq 0$
($m_\beta(\lambda-\sigma(\lambda)) \neq 0$) then $\lambda-\sigma(\lambda)$ is
either $-\alpha_i$ or $\theta,$ depending on which wall, and
$m_\theta(\lambda-\sigma(\lambda)) =1$
($m_\beta(\lambda-\sigma(\lambda)) =1$ if the root is short).  Thus
$N_{\lambda,\theta}^\lambda $ is  
at least $r$ minus the number of walls of the Weyl alcove to which  $\lambda$ is
adjacent  ($r_0$ minus
the number of walls dual to short roots to which $\lambda$ is
adjacent).   In the simply-laced case, only corners are adjacent to
$r$ walls.  In the nonsimply-laced case, only weights for which
$(\lambda,\alpha_i)=0$ for all short simple roots $\alpha_i$ are
adjacent to $r_0$ walls dual to a short root.
\end{proof}

\subsection{TQFTs from closed subsets}

\begin{lem}\label{lm:innerproduct}
If $k\lambda_i$ is invertible and $\gamma$ is any weight, then
$$S_{k\lambda_i,\gamma}/\qdim(\gamma)= C_{\phi_i(\gamma)} C_{k\lambda_i}^{-1}
C_\gamma^{-1}=e^{2 \pi i(\lambda_i,\gamma)}$$
where $\phi_i(\gamma)$ is the simple weight $k\lambda_i \trunc\gamma.$
\end{lem}

\begin{proof}
 First, we note that $\tau_i^{-1}(\lambda_i)=-\lambda_i,$ because they
have the same inner product with the simple roots, and that 
 $(\rho,\gamma-\tau_i(\gamma))=h(\lambda_i,\gamma)$ for
all $\gamma.$   It suffices to check the second claim on simple roots.   For $j\neq
i,$ we have that $\tau_i(\alpha_j)$ is another simple root of the same
length, different from
$\alpha_i,$ and thus both sides of the equation are zero (the object $\iota$).  
For
$\alpha_i,$ we have $\tau_i(\alpha_i)=-\theta,$ so
$(\rho,\alpha_i-\tau_i(\alpha_i))=(\rho,\theta)+1=h=h(\lambda_i,\alpha_i).$

By Equation (\ref{eq:theta})
 \begin{multline*}
C_{\phi_i(\gamma)} C_{k\lambda_i}^{-1}
C_\gamma^{-1}=\\
\exp\left(\frac{\pi i [(\tau_i(\gamma),\tau_i(\gamma)) -
(\gamma,\gamma) + 2(\rho,\gamma-\tau_i(\gamma))
-2k(\tau_i(\gamma),\lambda_i)]}{k+h} \right),
\end{multline*} 
where we have used the fact that
$\phi_i(\gamma)=k\lambda_i-\tau_i(\gamma).$  Using the fact that $\tau_i$ is an
isometry and the  identities in the previous paragraph gives
$$C_{\phi_i(\gamma)} C_{k\lambda_i}^{-1}
C_\gamma^{-1}=\exp\left(\frac{2\pi i}{k+h}[h(\lambda_i,\gamma)+
k(\lambda_i,\gamma)]\right)=e^{2 \pi i (\lambda_i,\gamma)}.\qed$$
\renewcommand{\qed}{\hbox{}}
\end{proof}

For the rest of this article let $\Gamma$ be a closed subset of the Weyl
alcove of the form $\Gamma_Z$ or $\Delta_Z.$

\begin{thm}\label{th:pseudoclass}
The degenerate
invertible elements of $\Gamma$ are exactly the invertible elements of $\Gamma$ which, viewed
under $e \circ k^{-1}$ as a subgroup of 
$(\Lambda/\Lambda_r)^*,$ annihilate the image of $\Gamma$ in $\Lambda/\Lambda_r.$
In particular if $\Gamma$ is $\Gamma_Z$ for some $Z,$ then the 
degenerate invertible elements are those 
in the image of $Z$  under $k\ell$ intersected with $\Gamma.$
These are even or odd depending on 
whether $k(\lambda_i,\lambda_i)$ is an even or odd integer.
\end{thm}

\begin{proof}
If $k\lambda_i$ is invertible, M\"uger shows it is degenerate for $\Gamma$ if and
only if $S_{k\lambda_i,\gamma}=\qdim(\gamma)$ for every $\gamma \in
\Gamma,$ which by   Lemma
\ref{lm:innerproduct}  is true if and only if
$e^{\lambda_i}$ annihilates all $\gamma \in \Gamma.$

To determine whether $k\lambda_i$ is odd or even, we need to check
whether $C_{k\lambda_i}= \pm 1,$ which is to say whether
$$k(k(\lambda_i,\lambda_i) + 2(\rho,\lambda_i))$$
is an even or odd multiple of $k+h.$  Notice $2(\rho,\lambda_i)=
(\rho,\lambda_i-\tau_i(\lambda_i)) = h (\lambda_i,\lambda_i)$ (Using
the identities in the proof of the previous lemma). So we
are asking whether $k(\lambda_i,\lambda_i)$ is an even or odd integer.
\end{proof}

\begin{cor}\label{cr:modular}
For every $\g$ and for every level $k,$ the ribbon category associated
to the full Weyl alcove is modular.
\end{cor}

\begin{cor}\label{cr:DW}
If $\Gamma=\Gamma_Z$  and $k$ is one of the levels conjectured by
Dijkgraaf and Witten to admit a Chern-Simons theory for $G,$ i.e., if
$k(\ell(z),\ell(z))/2$ is an integer for each $z$ in $Z,$
then in the ribbon category associated to $\Gamma$ all degenerate
objects are even invertible elements, and in fact form a group isomorphic to $Z.$
Thus as in Section 1 $\Gamma$ yields a modular category and a TQFT.
These levels are exactly 
the levels at which $Z$ embeds into $\Gamma_Z$ as even degenerate
objects via the map
$k\ell.$
\end{cor}

\section{Products of  Modular Categories and Tensor Products of
TQFTs}

At this point we have technically succeeded in the goal of the paper:
We have constructed  TQFTs associated to
nonsimply-connected groups at the levels predicted by physics.  But
in fact we have an embarrassment of riches, in that we have constructed
many more TQFTs than that.  First of all there are the closed subsets
$\Delta_Z$ of
the Weyl alcove consisting entirely of invertible elements, which at many values of
$k$ give TQFTs by the method of Section 1.  Second, the levels
$k$ suggested by Dijkgraaf and Witten are not the only ones giving
TQFTs by any means: While these authors suggest levels which are a
\emph{multiple} of a certain $N,$  it is easy to check from Theorem
\ref{th:pseudoclass} that $\Gamma_Z$ is modular (without the need for a
quotient) whenever $k$ is \emph{relatively prime} to $N$. More
generally, $\Gamma_Z$ gives a TQFT
whenever it contains no odd degenerate objects. It is
incumbent upon us to give as complete as possible a description of
these `unexpected' theories  and, if we claim to have verified the
expectations of physics, to show that they contain no new, nontrivial
information in some sense.   This is the goal of this section.

\subsection{Products of modular categories}

Suppose that $\Gamma$ is the label set of a ribbon category, and
$\Gamma$ contains two closed subsets $\Gamma'$ and $\Gamma''$ such that    
\begin{enumerate}
\item the intersection $\Gamma' \cap \Gamma''$ consists  of
even degenerate objects,
\item the product $\trunc$ of any element of $\Gamma'$ with an element
of $\Gamma''$ is simple  (i.e. is an element
of $\Gamma$),
\item every element of $\Gamma$ is a product of an element of
$\Gamma'$ and $\Gamma''$ and
\item if $\lambda' \in \Gamma'$ and $\lambda'' \in \Gamma''$ then
$C_{\lambda' \trunc \lambda''}=C_{\lambda'} C_{\lambda''}.$
\end{enumerate}
Then we say that $\Gamma$ is the \emph{product} of $\Gamma'$ and
$\Gamma''.$

Notice that $R_{\lambda',\lambda''}=R^{-1}_{\lambda'',\lambda'}$ because of
Condition 4.  In particular,  consider the invariant of a link with components
labeled by elements of $\Gamma.$   Every label can be written as a product of
a label in each of the subsets  (not uniquely, but since different
choices will disagree by a factor of an even degenerate object, it
will not effect the argument to follow), and using the tensor product property
of the link invariant, it can be written as the invariant of a link
with twice as many components, all labeled by elements of one of the
two subsets.   Now because of the condition on the $R$-matrix, this
invariant is equal to the invariant of two unlinked  copies of the
original link, one labeled by  labels in $\Gamma'$ and one by 
labels in $\Gamma''.$  In this sense the invariant associated to
$\Gamma$ is the product of the invariants associated to the two
factors.  This is the motivation for the definition.

The $S$-matrix for $\Gamma$ is the tensor product of the $S$-matrices
of the factors, so $\Gamma$ corresponds to a modular category if and
only if the factors do.

\begin{prop}\label{pr:product}
Suppose $Z \subset Z'$ are subgroups of the center of $G,$
 $\Gamma'$ is the closed subset
generated by $\Gamma_{Z'}$ and $\Delta_Z,$ and
$Z_0=Z \cap (k\ell)^{-1} [\Gamma_{Z'}].$ Then $\Gamma' \subset
\Gamma_{Z_0}$ is of the form $\Gamma_Y$ for some $Y$ and $\Delta_{Z_0}=\Delta_Z
\cap
\Gamma_{Z'}$ consists of degenerate invertibles for $\Gamma'.$  If all of
$\Delta_{Z_0}$ is even then 
$\Gamma'$ is the  product of $\Gamma_{Z'}$ and $\Delta_Z.$ These are
the only cases in which $\Gamma$ decomposes
into a product, apart from $D_{2n}.$  In the case $Z=Z',$ 
$\Gamma'=\Gamma_{Z_0}$ is one of the closed subsets described in Corollary \ref{cr:DW}. 
\end{prop}

\begin{proof}
Since $Z_0 \subset Z,$ the image of $Z_0$ under $\chi$  annihilates
$\Gamma_{Z}$ and hence $\Gamma_{Z'}.$  Since $Z_0 \subset (k
\ell)^{-1} [\Gamma_{Z'}],$
$\Delta_{Z_0} \subset \Gamma_{Z'}$ so the image of $Z$ under $\chi$ annihilates 
$\Delta_{Z_0},$ or equivalently the image of $Z_0$ annihilate
$\Delta_Z.$  Thus the 
image of $Z_0$ annihilates $\Gamma'.$ By Lemma \ref{lm:innerproduct},
this means $\Delta_{Z_0}$ (which of course is contained in $\Gamma'$)
will consist of degenerate units for $\Gamma'.$

If $\Delta_{Z_0}$ consists entirely of even degenerate objects, then
Condition 1 is met.  Conditions 2 and 3 are clear, and Condition 4
follows from Lemma \ref{lm:innerproduct}. In the case where $Z=Z',$
notice $\Gamma'$ contains $\Lambda_0 \cap \Lambda_r.$  Thus by Lemma
\ref{lm:exist} $\Gamma'$ consists of a union of cosets of $\Lambda_0/\Lambda_r,$
which is to say that it is of the form $\Gamma_{Y}$ for some
$Y.$  Of course $Y \supset Z_0,$ but anything in
$Y$ annihilates $\Gamma_Z$ and $\Delta_Z,$ so it is contained
in $Z$ and in $(k\ell)^{-1}[\Gamma_Z],$ so $Y=Z_0$ and
$\Gamma'=\Gamma_{Z_0}.$  
By Corollary 
\ref{cr:DW},  this is one of the Dijkgraaf-Witten theories.

To see these are the only cases of products of ribbon categories,
suppose $\Gamma$ is a closed subset which can be written as a product
of two closed subsets $\Gamma'$ and $\Gamma''.$  If both were of the
form $\Gamma_Z,$ and not of the form $\Delta_Z,$ 
then both would contain $\theta.$  Since this is not invertible except in
the case $\mathrm{su}_2$ at level $2,$ when one would have to be
$\Delta_{\Z_2},$ Condition 1 prevents both from being only of the form
$\Gamma_Z.$  They cannot both be of the form $\Delta_Z,$ because then
the product would be of the form $\Delta_Z$ with a $Z$ which was a
product of groups, which only happens for $\g=D_{2n}.$

Thus one must be of the form $\Delta_Z$ and one of the form
$\Gamma_{Z'}$ for some $Z, Z'.$  But by Condition 4 and Lemma
\ref{lm:innerproduct}, $Z$ must annihilate $\Gamma_{Z'}$ under
$\chi,$ and thus $Z \subset Z'.$ 
\end{proof}

\subsection{Tensor product of TQFTs}

\begin{prop}\label{pr:TQFTproduct}
If  the degenerate objects of $\Cat$ are all even and invertible and form a
cyclic group, and if the label set $\Gamma$ is a product of $\Gamma'$ and
$\Gamma'',$ then the modular quotient of $\Cat$ is the product of the
images in that quotient of $\Gamma'$ and $\Gamma''.$
 \end{prop}
\begin{proof}
A little thought will convince the reader that this result is almost
immediate assuming that the intersection (a group of degenerate
invertibles) acts freely on $\Gamma.$  This is in fact the case,
though one would like a more direct argument than the one below.
Let $Z_0$ be the group $\Gamma' \cap \Gamma''.$

Recall that, if $V$ is the vector space of formal linear combinations
of elements of $\Gamma,$ then one can extend the link invariant to an
invariant of links labeled by $V$ in such a way that it is linear in
each component.  Furthermore, the functor to $\Cat'$ gives a linear
map from $V$ into the corresponding $V'$ which is consistent with the
link invariant and a vector is in the kernel of this map if and only
if it is in the kernel of the invariant for every link. In particular
the Hopf link gives a nondegenerate pairing on the image of $V$ in
$V'.$ Brugui\`eres 
proves that if the Hopf link is labeled respectively by $\lambda \in
\Gamma$ and  $\omega= \sum_{\gamma \in \Gamma} \qdim(\gamma)
\gamma$  then the invariant is $0$ unless $\lambda$ is degenerate, in
which case it is $\qdim(\lambda) \qdim(\omega)=\qdim(\omega).$  The
same is necessarily true for $\lambda' \in \Gamma'$ and $\omega'=
\sum_{\gamma' \in \Gamma'} \qdim(\gamma') 
\gamma'$ and likewise for $\lambda'',\omega''$ in $\Gamma''.$  Thus the value
of  a Hopf link labeled by $\omega' \trunc \omega''$ and a typical element
$\lambda' \trunc \lambda''$ in $\Gamma$ is zero unless $\lambda'$ and
$\lambda''$ are both degenerate. In particular, since of course it is
nonzero if $\lambda' \trunc \lambda''$ is degenerate, every degenerate
object in $\Gamma$ can be written as a product of degenerate elements
in $\Gamma'$ and $\Gamma'',$ so in fact $\omega' \trunc \omega''$
gives zero on the Hopf link exactly when it is paired with a
nondegenerate simple object.  Thus a multiple of $\omega' \trunc
\omega''$ gives the same functional on $V$ via the Hopf link as
$\omega.$  Since the Hopf link labeled by $\omega $ and $\omega$ gives
$\qdim(\omega) |Z|$ and that labeled by $\omega$ and $\omega' \trunc
\omega''$ gives $\qdim(\omega) |Z|\cdot |Z_0|,$ we conclude that $\omega'
\trunc \omega''=|Z_0| \omega,$ and thus that $Z_0$ acts freely on
$\Gamma.$  

With this in hand, if $\lambda=\lambda' \trunc \lambda''$ then the
stabilizer of $\lambda$ is the product of the stabilizers of
$\lambda'$ and $\lambda''.$  The image of $\lambda'$ and $\lambda''$
in $\Cat'$ is a sum of as many simple objects as there are elements in
the stabilizer, and since each of these is a direct summand in the
image of $\lambda,$ they must each be simple.   Every simple
object in $\Cat'$ arises this way,  so  every simple object in
$\Cat'$ is the product of an object in the image of $\Gamma'$ and
$\Gamma''.$  Since $\Cat'$ contains no degenerate objects, it is
necessarily a product of these two subcategories.
\end{proof}

In Appendix \ref{A:product}, we define the notion of the tensor
product of two TQFTs and show that the product of two modular
categories gives the tensor product of their two TQFTs (that the
invariant is the product of the two invariants is an easy consequence
of the fact that the link invariant is the product of the link
invariants.    From this and the previous subsection we can deduce all
the tensor product relationships between the TQFTs we have constructed.
\begin{cor}
In the situation described in Proposition \ref{pr:product}, the three
TQFTs are related by $\ZZ_{\Gamma'}=\ZZ_{\Delta_Z} \tensor \ZZ_{\Gamma_{Z'}},$ and
(except for $\g=D_{2n}$) every TQFT arising from some $\Gamma_Z$ can be tensored
with one arising from 
$\Delta_{Z}$ to give one of the TQFTs conjectured by Dijkgraaf and
Witten.
\end{cor}

\begin{rem}
If we think of the product $\trunc$ as analogous to an algebra product, 
then to
say of two objects $\gamma$ and $\lambda$ that
$R_{\gamma,\lambda}=R_{\lambda,\gamma}^{-1}$ is analogous to saying that they
are commuting elements of the algebra.  Thus we should think of the subcategory
of degenerate objects as being the `center' of the ribbon category, and the
quotient of M\"uger and Brugui\`eres is the quotient of an algebra by its
center to give an algebra with a trivial center (modularity).   Extending this,
the definition of product is to say that the algebra is generated by two
commuting subalgebras, and thus the quotient is a product of their quotients.
\end{rem}

\subsection{Modular categories consisting entirely of invertibles}

It remains to identify the TQFTs $\ZZ_{\Delta_Z}$ where $Z$ is a subgroup of the
center $Z(G).$  In fact 
the resulting TQFT and three-manifold invariant have already been defined and
studied by Murakami, Ohtsuki and Okada \cite{MOO92}. 
\begin{prop}
If $\g \neq D_{2n}$ the three-manifold invariant arising from $\Delta_Z$ in the
case when
$\Delta_Z$ has no odd degenerate elements is
\begin{equation} \label{eq:moo}
Z_N(M,r)= \left( \frac{G_N(r)}{|G_N(r)|}\right)^{-\sigma(A)} |G_N(r)|^{-n}
\sum_{l \in (\Z_2)^n} r^{l^t A l}
\end{equation}
where $N$ is the order of $\Delta_Z$ modulo the even degenerate objects, 
$A$ is the linking matrix of a framed link presenting the three-manifold
$M,$
$n$ is the number of components, $\sigma(A)$ is the signature of $A,$ $l^t$ is
the transpose of the vector $l$ and $G_N(r)= \sum_{m=1}^N r^{m^2}$ where $r$ is
$\exp(k \pi i (\lambda_i,\lambda_i)),$ with $\lambda_i$ the element of $\ell[Z]$
giving $k(\lambda_i,\lambda_i)$ the largest denominator.  This is the
invariant constructed by Murakami, Ohtsuki and Okada.  
\end{prop}

\begin{proof}
 
Of course it suffices to see that, with $I$ as defined in Section
1.3,  $I(L)$ is $\sum_{l \in (\Z/N)^n}
r^{l^t A l}.$  

Let $\lambda$ be a generator of  $\Delta_Z.$  
then  $\lambda^N$ is even
degenerate and no smaller power of $\lambda$ is degenerate. 
 To say that
$\lambda^n$ for some $n$ is degenerate for $\Delta_Z$ is to say that
$C_{\lambda^{n+1}} =C_{\lambda^n} C_\lambda,$   since it suffices to check
the degeneracy condition against a generator.   Now
$\lambda=k\lambda_i$ and $\lambda^n=k\lambda_j$ for some $\lambda_i, \lambda_j
\in
\ell[Z],$ so by Lemma \ref{lm:innerproduct}, this is to say that
$k(\lambda_i,\lambda_j)$ is an integer, which by the fact that the map $e$ is a
homomorphism is equivalent to saying $kn(\lambda_i,\lambda_i)$ is an integer.  Thus
for any $\lambda_j$ the denominator of $k(\lambda_j,\lambda_j)$ represents the order
of $k\lambda_j$ in the quotient group $\Delta_Z$ modulo even
degenerates, and thus the statement that $\lambda$ is a generator  
is equivalent to $k(\lambda_i,\lambda_i)$ having maximal denominator
as in the statement of the
proposition.   Let $r=\exp(\pi i k(\lambda_i,\lambda_i)).$  Notice $r^{2N}=1, $ and
$N$ is the least natural number for which this is true.

The fact that $\lambda^N$ is even degenerate means that $C_{\lambda^N}=1.$ 
Applying Lemma \ref{lm:innerproduct} recursively shows $C_{\lambda^N}=C_\lambda^N
\exp( \pi i k N (N-1) (\lambda_i,\lambda_i)),$ while of course
$C_\lambda=\exp(\pi i (k\lambda_i,k\lambda_i+2\rho)/(k + h))=\exp(\pi i k
(\lambda_i,\lambda_i))=r,$ so $C_{\lambda^N}=r^{N^2}.$   In order for this to be $1$
we conclude that $N$ is odd and $r$ is a primitive $N$th root of unity or $N$ is even
and $r$ is a primitive $2N$th root of unity.   We claim that the link
invariant of a link with linking matrix $A$ and with the $n$
components labeled by the vector of labels $l \in (\Z/N)^n$ or $l
\in (\Z_{2N})^n$ depending on the parity of $N,$ where we mean
labeling the $i$th component by $\lambda^{l_i},$ is $r^{l^t A l}.$   

To see this, notice that because all labels are units $R_{\lambda^n,\lambda^m}
R_{\lambda^m,\lambda^n}=S_{\lambda^n,\lambda^m} \cdot \mathrm{id}= r^{2nm},$ which means
that switching an undercrossing to an overcrossing in a  link projection with these
labels multiplies the link invariant by $r^{2nm},$  which is exactly the effect this
move has on $r^{l^t A l}.$   Since such moves will untie any link to
a collection of unlinked framed unknots, it suffices to check that the
formula is correct
 on these.   This follows from the fact that both assign $r^{mn^2}$ to the
$m$-framed unknot with label $n.$  

Now that we know the link invariant, $I(L)$  is either  $\sum_{l \in
(\Z_N)^n} r^{l^t A l}$ or $\sum_{l \in (\Z_{2N})^n} r^{l^t A l}$
depending on the parity of $N.$   In the second
case the $\Z_2$ symmetry of the labels because of the even
degenerate unit means that this sum is $2^n  \sum_{l \in
(\Z_N)^n} r^{l^t A l}$ and as noted in Section 1.3 an overall
factor in $I$ depending on the number of components is canceled out
in the three-manifold invariant.
\end{proof}

\begin{rem}
The product of this invariant with its conjugate (or on the level of TQFTs, the
tensor product of this TQFT with its conjugate TQFT) is an  example of the
finite group invariants constructed by Dijkgraaf and Witten in the same paper
\cite{DW90} as the one where they construct  Chern-Simons theory from
nonsimply-connected Lie groups.    Thus Murakami et al's theories stand in the
same relationship to the finite group theories as Chern-Simons theory does to
Turaev-Viro  (sometimes called three-dimensional gravity).
\end{rem}
\section*{Appendix}

\renewcommand{\thesection}{A}
\setcounter{subsection}{0}

\subsection{Constructing TQFTs from the ribbon $*$-category}\label{A:TQFT}

Recall the definition of a $2$-framed three-dimensional  TQFT given in
\cite{Sawin99}, an assignment 
of a vector space
$\ZZ(\Sigma_g)$ to a given oriented $2$-framed surface $\Sigma_g,$ and
an assignment of a functional $\ZZ(M)\colon\bigotimes_{i=1}^n \ZZ(\Sigma_{g_i})
\to \C$ to each $2$-framed three-manifold with boundary parameterized
by an orientation reversing homeomorphism from $\bigcup_{i=1}^n
\Sigma_{g_i}$ satisfying certain conditions. 

Following \cite{Sawin96a}  we construct a TQFT from a modular category
as follows.  Choose a handlebody $H_g$ with boundary $\Sigma_g,$
choose a framed graph in $H_g$ whose image is a retract of $H_g,$ and
let $\ZZ(\Sigma_g)$ be the set of all labelings of the graph, with two
labelings identified if they give the same invariant for all
embeddings  (this is just the sum over all labelings of the edges of
the tensor product over all vertices of the space of possible
labelings of the vertex, as described in Section 1.1).  If $M$ is a
$2$-framed three-manifold with boundary parameterized by $\bigcup_{i=1}^n
\Sigma_{g_i},$ we can present $M$ by an embedding of $\bigcup_{i=1}^n
H_{g_i}$ into $S^3$ together with a framed link in the complement of
the embedding, with $M$ homeomorphic to surgery on the link in the
complement of the embedding by a map preserving the parameterization
of the boundary.   The functional $\ZZ(M)$ evaluated on vectors
corresponding to certain labelings of the $n$ graphs in $\{H_{g_i}\}$
is the invariant of the labeled graphs embedded in $S^3$ together
with the link, each component of the link being labeled by $\omega$
(times a normalization factor which need not concern us).  

We will identify a labeling of a graph by data in the modular quotient
$\Cat'$ 
with a labeling of a slightly modified graph by data in the original
category $\Cat$ in such a way that the invariant of any ribbon graph labeled
by data in the modular quotient is the same as the invariant of the
modified graph with the associated data in the original graph
(actually, it will be a formal linear combination of such labelings).
This 
will tell us how to construct the vector spaces $\ZZ(\Sigma_g)$ and
the functionals $\ZZ(M)$  from data in $\Cat$  (we have already seen
that links labeled by $\omega$ in $\Cat'$ can be replaced by the same
link labeled by $\omega$ in $\Cat$ up to a factor, and that continues
to be true in the presence of other graph components).

Given a graph $G,$ the modified graph will be formed by adding  a new edge to
each vertex of $G,$
 with all the new edges meeting in a single new vertex.  Of
course an embedding of $G$ into $S^3$ extends in many ways to an
embedding of this larger graph, but since
 all the new edges will be labeled by degenerate objects in
$\Cat,$  all these embedding will yield the same invariant and we are
not obliged to specify a particular embedding.  

Recall that M\"uger's construction of $\Cat'$ proceeds in two steps.  He
first constructs an intermediate category $\Cat^0$ whose objects are
the objects of $\Cat$ but whose morphisms from $\lambda$ to $\gamma$
are $\bigoplus_\mu \hom(\lambda,\gamma \trunc \mu),$ where the sum is
over all degenerate simple objects $\mu\in \Gamma.$   The category $\Cat'$ is
constructed from this by the usual process of closing under
subobjects.   In particular $\Cat^0$ is a full subcategory of $\Cat'$
and for each  object $\lambda$ of $\Cat'$ there
is an object $\gamma$ of $\Cat^0$ and a morphism $f\colon\lambda \to
\gamma$ such that $f^*f$ is the identity and $ff^*$ is a minimal
idempotent.  

In particular, if $G$ is a graph labeled by data in $\Cat'$ we can
replace each object $\lambda$ labeling an edge by the corresponding
$\gamma$ in $\Cat^0$ as above, and compose each vertex label with the
appropriate $f$ and $f^*$ morphisms to get a new labeling by data in
$\Cat^0$ with the same value of the invariant on each embedding.  Thus
we need only describe the process for a graph labeled by data in
$\Cat^0.$  

So let $G$ be an abstract graph labeled by data in $\Cat^0,$ and let
$G'$ be the extended version of $G$ above.   Each edge of $G$ is
labeled by an object of $\Cat^0,$ which is also an object of $\Cat,$
so we label the corresponding edge of $G'$ by that object.   Assume
that each new edge in $G'$ was added so as to be last in the ordering
of edges around its vertex.   Then the label of a given vertex is an
element of $\bigoplus_\mu \hom(\lambda \trunc \mu^\dagger,\iota) \iso
\bigoplus_\mu\hom(\lambda,\mu),$ where $\lambda$ is the 
tensor product of the edge labels of $G$ or their duals as in Section 1.1.
We can assume that each label is in exactly one of these direct
summands, as the general case can be expressed as a linear combination
of such.  So for each vertex there is a degenerate label $\mu$ such
that the the label is  $x\in\hom(\lambda,\mu).$  Label the new edge of
$G'$ by $\mu$ (oriented away) 
and the vertex by $x.$   The new  vertex needs a label in
$\hom(\mu_1 \trunc \cdots \trunc \mu_n,\iota),$  which of course is the identity
morphism if $\mu_1 \trunc \cdots \trunc \mu_n=\iota$ and $0$ otherwise.  It is
now a simple exercise to check that for any embedding of $G$ and any
extension of that to an embedding of $G',$ the invariant of $G$ equals
the invariant of $G'.$

Two other observations follow by simple calculations.  The first is
that we could as well choose $G'$ to have only one new edge for each
connected component, the cost being the labels on the edges of $G$
might have to be multiplied by degenerate objects.  Thus the vector
space $\ZZ(\Sigma_g)$ is still spanned by labelings of a graph in
$H_g$ which is a retract of $H_g,$ but it has an extra edge (which is
required to be labeled by a degenerate object) and an
extra univalent vertex.  

 In particular the vector space of the torus, which can be viewed as
the space of labels for the link invariant, will have a basis element
for each orbit under $Z$  of simple objects on $\Cat,$ as well as a
basis element for each such object and any nontrivial element of its
stabilizer, the stabilizing element labeling the extra edge.   These
basis elements with nontrivial labels on the extra edge give a nonzero
link invariant only when there is another component in the link which
is labeled by a basis element with a nontrivial label.  A similar
picture appears in the literature in descriptions of WZW models
arising from nonsimply-connected Lie groups  \cite{MS89,GW86}.  

Note that each of the vector spaces associated to surfaces decomposes
naturally as a sum of sectors according to the label of the extra
edge.  This label is an element of the group $Z,$ which in the
nonsimply-connected groups we explore in Section 2 represents the
fundamental group of the group, which also indexes different principal
bundles of the group $G_Z$ over a connected surface.  Thus we may view the
vector spaces $\ZZ(\Sigma_g)$ as a sum over contributions from the
different principal bundles of the surface, as expected from the
physics.

The second observation is that if the process of constructing $G'$
with labels in $\Cat$ is applied to a link component labeled by
$\omega$ (think of this as one edge labeled by the direct sum of all
simple representations and a single bivalent vertex labeled by the sum
of the canonical duality map times the quantum dimension), the
labeling morphisms lie in the trivial component of
$\bigoplus_\mu\hom(\lambda,\mu),$ and thus we do not need to add the
extra edge, and the result is the same link component labeled by
$\omega$ is $\Cat.$   Thus although $\Cat$ is not modular $\omega$
continues to be a perfectly good `surgery label,' and the only
adjustment that needs to be made to Reshetikhin and Turaev's process is
an extension of the graphs spanning $\ZZ(\Sigma_g).$

\subsection{The product of two modular categories gives the tensor
product of their TQFTs}
\label{A:product}

 We say that a TQFT $\ZZ$
is the  tensor product of
two TQFTs $\ZZ_1$ and $\ZZ_2$ if for each $g$ there is an isomorphism
$\Phi_g\colon\ZZ(\Sigma_g) \to \ZZ_1(\Sigma_g) \tensor \ZZ_2(\Sigma_g)$ such
that $\ZZ(M)\circ \Phi^{-1}=\ZZ_1(M) \tensor \ZZ_2(M),$ where here we
have identified the two domain spaces which are connected by a
rearrangement of the tensor factors.

\begin{prop}
If $\Cat$ is modular and has label set $\Gamma$ which is a product of
$\Gamma'$ and $\Gamma'',$ then the TQFT associated to $\Cat$ is the
tensor product of the TQFTs associated to the subcategories determined
by $\Gamma'$ and $\Gamma''.$
\end{prop}

\begin{figure}[hbt]
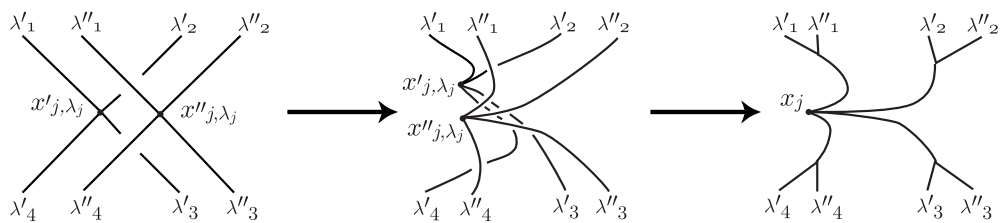


\caption{Pictorial version of the $I$ isomorphism}\label{fg:transform}
\end{figure}

\begin{proof}
Of course if $\Cat$ is modular then $\Gamma' \cap \Gamma'' = \{\iota\}$ and
it is easy to see that the category is a product of the two
corresponding subcategories.  

The central observation is that if $G$ is any abstract  graph
and $l$ is a choice of an object $\gamma_i \in \Gamma$ for each edge
$i$ and $x_j \in \hom(\lambda_j,\iota)$ for each vertex $j,$ where
$\lambda_j$ is constructed out of $\{\gamma_i\}$ as in the definition
of the graph invariant in Section 1.1, then there exists a set of pairs
of labels $l_\beta'$ and $l_\beta''$ for some index $\beta,$ with all
the data in $l'$ and $l''$ coming from the categories associated to
$\Gamma'$ and $\Gamma''$ respectively, such that for any embedding of
$G$ into $S^3$ the invariant of $G$ labeled by $l$ is the sum over all
$\beta$ of the product of the invariants of the embedding of $G$
labeled by $l_\beta'$ and $l_\beta''$ respectively.

To see this, notice there are unique $\gamma_i'\in \Gamma'$ and
$\gamma_i''\in \Gamma''$ such that $\gamma_i=\gamma_i' \trunc
\gamma_i''.$  If we write $\lambda_j= \overline{\gamma}_{i_1} \trunc
\overline{\gamma}_{i_2} \trunc \cdots \trunc \overline{\gamma}_{i_n},$
with $\overline{\gamma}_i=\gamma_i$ or $\gamma_i^\dagger,$ and define
$\lambda_j' =\overline{\gamma}'_{i_1} \trunc \cdots \trunc
\overline{\gamma}'_{i_n}$ and $\lambda_j''=\overline{\gamma}''_{i_1}
\trunc \cdots \trunc \overline{\gamma}''_{i_n},$ then there is a
canonical isomorphism formed from a product of $R$-morphisms
$\lambda_j' \trunc \lambda_j'' \iso \lambda_j$ because of the relation
$R_{\gamma',\gamma''}=R_{\gamma'',\gamma'}^{-1}.$  In particular there
is an isomorphism
$$I\colon\hom(\lambda_j',\iota) \tensor \hom(\lambda_j'',\iota) \to
\hom(\lambda_j,\iota)$$
so we can write 
$$x_j= I(\sum_{\alpha_j} x_{j,\alpha_j}' \tensor x_{j,\alpha_j}'').$$
Now for a choice $\beta$ of an $\alpha_j$ for each vertex $j,$ we can
define $l_\beta'$ and $l_\beta''$ by $\{\gamma_i',x_{j,\alpha_j}'\}$
and $\{\gamma_i'',x_{j,\alpha_j}''\}$ respectively.  Consider an
embedding of $G,$ and consider two parallel copies of that embedding
of $G$  one shifted from the other by the framing, and one labeled by
$l_\beta',$ the other by $l_\beta''.$  On the one hand the fact that
$R_{\gamma',\gamma''}=R_{\gamma'',\gamma'}^{-1}$ means we can pass
these two copies of $G$ through each other to obtain two entirely
disjoint copies of $G$ with the same invariant, and thus the invariant
of the doubled graph is the product of the invariant of the two
labeled graphs separately.  On the other hand, since the
transformation of Figure \ref{fg:transform} does not change the
invariant, it equals the invariant of one copy of the embedded $G$
with edges labeled by $\gamma_i$ and vertices labeled by
$I(x_{j,\alpha_j}' \tensor x_{j,\alpha_j}'').$  Summing over all
$\beta$ we get the invariant of $G$ labeled by $l.$

Now $\ZZ(\Sigma_g)$ is the space spanned by all nonzero labelings of a
fixed graph embedded in a handlebody $H_g$ with boundary $\Sigma_g$
whose image is a retract of the handlebody.  The map $I$ defined
above then gives a linear map $\Phi\colon\ZZ(\Sigma_g) \to \ZZ'(\Sigma_g)
\tensor \ZZ''(\Sigma_g)$ which is easily seen to be an isomorphism
(the map $I$ above provides the inverse).

If $M$ is a $2$-framed manifold presented by an embedding of
$\bigcup_{i=1}^n H_{g_i}$ into $S^3$ and a surgery link on the
complement, then the value of the functional on a vector given by
labelings of graphs in $H_{g_i}$ is (up to some normalizations) the
invariant applied to the 
embeddings of these graphs together with the surgery link, with each
component labeled by $\omega.$  Now notice that since $\Gamma=\Gamma'
\times \Gamma''$ we have $\omega=\omega'\trunc \omega''=\Phi(\omega'
\tensor \omega'').$  Thus the image under $\Phi$ of the tensor product
of vectors coming from labelings which compute $\ZZ'(M)(v)$ and
$\ZZ''(M)(w)$ is the vector coming from the labeling which computes
$\ZZ(M)(\Phi(v \tensor w)),$ so
$$\ZZ(M)\circ \Phi= \ZZ'(M) \tensor \ZZ''(M)$$
after the appropriate identification of the domain spaces.
\end{proof}
$$\pic{transform}{80}{-30}{0}{0}$$


\end{document}